\documentstyle[11pt]{article}
\normalbaselines

\newtheorem{thm}{Theorem}[section]
\newtheorem{lem}{Lemma}[section]
\newtheorem{prop}{Proposition}[section]
\newtheorem{cor}{Corollary}[section]

\begin{document}

\def\bea*{\begin{eqnarray*}}
\def\eea*{\end{eqnarray*}}
\def\ba{\begin{array}}
\def\ea{\end{array}}
\count1=1
\def\be{\ifnum \count1=0 $$ \else \begin{equation}\fi}
\def\ee{\ifnum\count1=0 $$ \else \end{equation}\fi}
\def\ele(#1){\ifnum\count1=0 \eqno({\bf #1}) $$ \else \label{#1}\end{equation}\fi}
\def\req(#1){\ifnum\count1=0 {\bf #1}\else \ref{#1}\fi}
\def\bea(#1){\ifnum \count1=0   $$ \begin{array}{#1}
\else \begin{equation} \begin{array}{#1} \fi}
\def\eea{\ifnum \count1=0 \end{array} $$
\else  \end{array}\end{equation}\fi}
\def\elea(#1){\ifnum \count1=0 \end{array}\label{#1}\eqno({\bf #1}) $$
\else\end{array}\label{#1}\end{equation}\fi}
\def\cit(#1){
\ifnum\count1=0 {\bf #1} \cite{#1} \else 
\cite{#1}\fi}
\def\bibit(#1){\ifnum\count1=0 \bibitem{#1} [#1    ] \else \bibitem{#1}\fi}
\def\ds{\displaystyle}
\def\hb{\hfill\break}
\def\comment#1{\hb {***** {\em #1} *****}\hb }

\newcommand{\TZ}{\hbox{\bf T}}
\newcommand{\MZ}{\hbox{\bf M}}
\newcommand{\ZZ}{\hbox{\bf Z}}
\newcommand{\NZ}{\hbox{\bf N}}
\newcommand{\RZ}{\hbox{\bf R}}
\newcommand{\CZ}{\,\hbox{\bf C}}
\newcommand{\PZ}{\hbox{\bf P}}
\newcommand{\QZ}{\hbox{\bf Q}}
\newcommand{\HZ}{\hbox{\bf H}}
\newcommand{\EZ}{\hbox{\bf E}}
\newcommand{\GZ}{\,\hbox{\bf G}}

\font\germ=eufm10
\def\goth#1{\hbox{\germ #1}}
\vbox{\vspace{38mm}}

\begin{center}
{\LARGE \bf Fermat Surface and Group Theory in Symmetry of Rapidity Family in Chiral Potts Model} \\[10 mm] 
Shi-shyr Roan \\
{\it Institute of Mathematics \\
Academia Sinica \\  Taipei , Taiwan \\
(email: maroan@gate.sinica.edu.tw ) } \\[25mm]
\end{center}

\begin{abstract}
The present paper discusses various mathematical aspects about the rapidity symmetry in chiral Potts model (CPM) in the context of algebraic geometry and group theory . We re-analyze the symmetry group of a rapidity curve in $N$-state CPM, explore the universal group structure for all $N$, and further enlarge it to modular symmetries of the complete rapidity family in CPM. As will be shown in the article that all rapidity curves in $N$-state CPM constitute a Fermat hypersurface in $\PZ^3$ of degree $2N$ as the natural generalization of the Fermat K3 elliptic surface $(N=2)$, we conduct a thorough algebraic geometry study about the rapidity fibration of Fermat surface and its reduced hyperelliptic fibration via techniques in algebraic surface theory. Symmetries of rapidity family in CPM and hyperelliptic family in $\tau^{(2)}$-model are exhibited through the geometrical representation of the universal structural group in mathematics.

\end{abstract}
\par \vspace{5mm} \noindent
  \par \noindent
{\rm 2008 PACS}: 05.50.+q, 02.20.Bb, 02.40.Tt \par \noindent

\setcounter{section}{0}
\section{Introduction}
\setcounter{equation}{0}
In the study of the two-dimensional $N$-state chiral Potts model (for a brief
history account, see, e.g. \cite{B91} or \cite{Mc} section 4.1 and references therein), 
"rapidities"  of the statistical model are 4-vector ratios $p=[a, b, c, d]$ in the projective 3-space $\PZ^3$ satisfying the following equivalent equations of an algebraic curve of genus $N^3-2N^2+1$:
\bea(llll)
{\goth W}_{k', k}  : &    \left\{ \begin{array}{l}
ka^N + k'c^N = d^N , \\
kb^N + k'd^N = c^N . \end{array}  \right.  & \Longleftrightarrow &   
\left\{ \begin{array}{l}
a^N + k'b^N = k d^N  , \\
k'a^N + b^N =
k c^N , \end{array}  \right.
\elea(Wk'k)
where $k', k$ are temperature-like parameters $\neq 0, \pm 1$ with $k^2 + k'^2 = 1$. This high genus rapidity curve possesses a large symmetry group, which has played an essential role in the solvability of the $N$-state chiral Potts model (CPM) as the natural generalization of Ising-model ($N=2$). It is easy to see that the following transformations of $\PZ^3$ preserve the equation (\req(Wk'k)):
\bea(lll)
M^{(1)}: [a,b, c, d] \mapsto [\omega a, b, c, \omega d ] , & 
M^{(2)}: [a,b, c, d] \mapsto [\omega a, \omega b, c,  d ] , \\
M^{(3)}: [a,b, c, d] \mapsto [ c, \omega^{\frac{1}{2}}d, \omega^{\frac{-1}{2}}a, \omega^{-1}b ] , &
M^{(4)}: [a,b, c, d] \mapsto [ a, b, \omega^{-1}c, d ] , \\
M^{(5)}: [a,b, c, d] \mapsto [ d, \omega^{\frac{1}{2}}c, \omega^{\frac{-1}{2}}b, a ], & M^{(0)}: [a,b, c, d] \mapsto [\omega a, b, \omega c, d ], 
\elea(Sym)
where $\omega = e^\frac{2 \pi {\rm i}}{N}$, hence induce the symmetry of ${\goth W}_{k', k}$ as elements in the automorphism group ${\rm Aut}({\goth W}_{k', k})$. Indeed, for $N \geq 3$, ${\rm Aut}({\goth W}_{k', k})$ is an order $4N^3$ group generated by automorphisms in (\req(Sym)) (\cite{AuP}, \cite{R04}  Proposition 1):
\bea(lll)
{\rm Aut}({\goth W}_{k', k})= \langle M^{(j)}\rangle_{j=0}^5,  & |{\rm Aut}({\goth W}_{k', k})| =4N^3, & (N \geq 3).
\elea(AutWk)
The purpose of this paper is to explore the universal structure of ${\rm Aut}({\goth W}_{k', k})$ for all $N$, based in part on earlier quantitative studies in \cite{R04}, then enlarge the ${\goth W}_{k', k}$-symmetry of a single curve to "modular" symmetries of the complete family of rapidity curves for all $(k', k)$ including the degenerated $k'=0, \pm 1, \infty$. 
The parameters $(k', k)$ of (\req(Wk'k)) can be identified with the 3-vector ratios $[k', k, 1] \in \PZ^2$ in a quadratic hypersurface of $\PZ^2$, which is biregular to $\PZ^1$ via the correspondence
\bea(llll)
 \Lambda := \{ (k', k) \in  \PZ^2 | k^2 + k'^2 = 1 \} & \simeq & \{ \kappa' \in \PZ^1 \} , &  
\kappa' = k' + {\rm i} k = (k' - {\rm i} k)^{-1} = (\frac{k' + {\rm i} k}{k' - {\rm i} k})^\frac{1}{2}.
\elea(kk'kap)
Here $\kappa'$  is a complex parameter including $\infty$, identified with elements in $\PZ^1$ via
$$
\CZ \ni \kappa'  \longleftrightarrow [\kappa, 1] \in \PZ^1, ~ ~ \kappa'= \infty  \longleftrightarrow [1, 0] \in \PZ^1.
$$
Note that in (\req(kk'kap)), $k'= \infty$ iff $k= \infty$ with $\frac{{\rm i} k}{k'}= 1$ or $-1$, which corresponds to $\kappa' = \infty, 0$ respectively. Then for each $(k', k) \in  \PZ^2$ in (\req(kk'kap)), the four equations in (\req(Wk'k)) defines a rapidity curve in $\PZ^3$, i.e. other than those in (\req(Wk'k)) for $k' \neq 0, \pm 1, \infty$, the rest are degenerated curves consisting of $N^2$ lines of $\PZ^3$ defined by 
\bea(llll)
{\goth W}_{k', k}, &  (k', k) =(0, \pm 1), \kappa'= \pm {\rm i} & :  
a^N  = \pm d^N , &   
b^N = \pm c^N ; \\
{\goth W}_{k', k},  & (k', k) =(\pm 1, 0), \kappa'= \pm 1 &:   
 c^N = \pm  d^N , &  a^N = \mp b^N  ; \\
{\goth W}_{k', k}  & (k', k) =(\infty,  \infty_\pm ), \kappa'= \infty^{\pm 1}&:    a^N  = \mp  {\rm i} c^N , &   b^N = \mp {\rm i}  d^N ,
\elea(Wdeg)
where $(k', k) =(\infty,  \infty_\pm )$ denotes the case $(k',  \frac{{\rm i}k}{k'})\rightarrow (\infty,  \pm 1)$ as $k' \rightarrow \infty$. The collection of all rapidity curves in (\req(Wk'k)) and (\req(Wdeg)) form a family over $\Lambda$:
\bea(lll)
{\goth W} ( = {\goth W}^{(N)} ) = \bigsqcup_{(k' k) \in \Lambda } {\goth W}_{k', k}  &\longrightarrow &\Lambda ,
\elea(WFaml)
which can be identified with a Fermat hypersurface in $\PZ^3$ of degree $2N$ (see, (\req(WFerm)) in the paper). Hence the symmetry of Fermat rapidity fibration (\req(WFaml)) over $\Lambda$ can be studied in the context of algebraic geometry as the natural generalization of the elliptic fibration of Fermat K3 surface ($N=2$). In this work, we identify the symmetry group of the rapidity family (\req(WFaml)) as an extension of the automorphism group (\req(AutWk)) with $PSL_2(\ZZ_4)$ generated by two "modular" symmetries $T, S$ of the fibration (\req(WFaml)) (see, Theorem \ref{thm:WGtN} in the paper). 
Due to the lack of difference property of rapidities in (\req(Wk'k)) for $N \geq 3$, the computation of some physical  interesting quantities in CPM, such as eigenvalue spectrum and eigenvectors \cite{AMP, B90, B93, MR, R905, R1003} or order parameter \cite{B05a, B05}, relies on the functional-relation method in \cite{BBP}, by regarding CPM as a descendant of the (six-vertex) $\tau^{(2)}$-model \cite{BazS} (see, also \cite{R1206}), with the $\tau^{(2)}$-spectral parameter lying on some hyperelliptic curves reduced from (\req(Wk'k)). Through the principle of symmetry reduction, the algebraic geometry study of the rapidity fibration (\req(WFaml)) can be carried over to the complete family of hyperelliptic curves. We are able to determine the geometrical properties and symmetries of the hyperelliptic fibration, especially the singularity structure around degenerated curves, by techniques of surface theory in algebraic geometry. In this paper, we conduct a qualitative investigation of symmetry about the chiral Potts rapidity family in (\req(WFaml)), in contrast to the usual quantitative approach, where the discussion proceeds through the explicit form of symmetries as in (\req(Sym)). The common structural characters of the symmetry group for all $N$ will be our main concern. Indeed, we build a mathematical model in group theory, which contains the universal structure of automorphism groups of the chiral Potts curve in (\req(Wk'k)) and the rapidity fibration in (\req(WFaml)). Through the representation of the universal group, one obtains the quantitative expression of symmetries of the rapidity family (\req(WFaml)), as well reproduces accurately the known ones in (\req(Sym)).

This paper is organized as follows. In Section \ref{sec:Alg}, we setup a mathematical model in the context of group-theory formulation for the study of symmetries of chiral Potts curves and rapidity family.
The solvable groups and its $PSL_2(\ZZ_4)$-extension, ${\goth G} \subset \widetilde{\goth G}$, ${\goth G}_N \subset \widetilde{\goth G}_N ~ (N \geq 2)$ (see, (\req(CPG)) (\req(uCPm)) (\req(CPNg)) (\req(CPNmg)) 
in the paper), are introduced in Subsection \ref{ssect:UCPG},  where the structures and properties of the groups are discussed in detail by using the approach of mathematical derivations.
In Subsection \ref{ssec:Dih}, we discuss the $(\ZZ_2 \times D_N)$-structures as quotients of the CP group ${\goth G}_N$ in group theory, where $D_N$ is the dihedral group which, together with $\ZZ_2$, describes the symmetries of  hyperelliptic curves in $\tau^{(2)}$-model of CPM. The relationship of these $(\ZZ_2 \times D_N)$-structures under the action of modular symmetries is also examined through the group structure of $\widetilde{\goth G}_N$.
Section \ref{sec.Fermat} is devoted to the algebraic geometry study of rapidity family (\req(AutWk)) in CPM and the hyperelliptic family in $\tau^{(2)}$-model. In Subsection \ref{subs.FSym}, we first precisely identify the rapidity family (\req(AutWk)) of $N$-state CPM with a degree-$2N$ Fermat hypersurface of $\PZ^3$ (see, (\req(WFerm)) in the article). Through the representation theory of  structure groups introduced in Section \ref{sec:Alg}, the automorphism group of ${\goth W}_{k', k}$ in (\req(Wk'k)) is correctly reproduced via the geometrical representation of ${\goth G}_N$, as well as the modular symmetries of (\req(AutWk)) which generate the extended group $PSL_2(\ZZ_4)$ in $\widetilde{\goth G}_N$. Furthermore, $\widetilde{\goth G}_N$ is geometrically characterized as the automorphism group of the rapidity-fibered Fermat surface (Theorem \ref{thm:WGtN} in the paper). In Subsection \ref{subs.hypell}, we investigate the geometry and symmetry structure of the hyperelliptic(-curve) families in CPM, which are related to the rapidity fibration of Fermat surface  by the symmetry reduction. However the global structure of the hyperelliptic family is drastically affected by the reduction process, in particular, the orbifold singularities that occurs in the degenerated fibers. By using techniques in algebraic surface theory and toric geometry, we explicitly construct the minimal resolution of the hyperelliptic-fibered surface. An analysis about the geometry properties and global symmetries of the surface is performed in details by the method of algebraic geometry.  In Subsection \ref{subs.K3}, we focus on the $N=2$ case, where the rapidity family is the well-known Fermat elliptic K3 surface, and the rapidity discussions in previous sections (in Ising model case) can be also illustrated in the context of uniformalization of elliptic curves. The demonstration provides a conceptual insight about the group-theory approach of CPM symmetry in this work that deserves to be known, as well as some additional informations valid only for $N=2$. By using the theta-function representation of rapidity variables, the symmetry of  the Fermat K3 elliptic fibration is well described by elliptic and modular transformations of uniformalization parameters. Geometrically, all rapidity families for $N=2$ are elliptic K3 surfaces over $\HZ/PSL_2(\ZZ_4)$.

Notation:  In this paper, we use the standard notations in group theory. Let $G$ be a group, and $H, H'$ be  subgroups, $B$ be a subset of $G$. We denote
\bea(lll)
\langle B \rangle &:= {\rm ~ the ~ subgroup ~ of ~} G ~ {\rm ~ generated ~ by ~} B, \\
\langle \langle B \rangle \rangle_{\rm N}&:= {\rm the ~ normal ~ subgroup ~ of ~} G ~ {\rm ~ generated ~ by ~} B , \\
N(H, H') &:= \{ g \in G | g H g^{-1} = H' \}, \\
N(H) &:= N(H, H) ~ ~{\rm ~ the ~ normalizer ~ of ~} H ~ {\rm in ~ } G.
\elea(BNH)

\section{Algebraic Theory in Symmetry of Chiral Potts Model \label{sec:Alg}}
\setcounter{equation}{0}
In this section, we build a mathematical model in group theory, based on common properties of rapidity symmetries in CPM.  The algebraic formulation of the symmetry groups contains the universal structure and essential characters of automorphism groups appeared in the study of CPM. However, in contrast to the usual quantitative approach in CPM, discussions in this section are carried out in a form of abstract mathematical derivation in group theory, knowledges in rapidities and CPM not required. The connection between the abstract groups in mathematic and the symmetry in CPM will be discussed later in Section \ref{sec.Fermat}.

\subsection{ CP group and modular CP group \label{ssect:UCPG}}
First we define the group which characterizes the universal structure of automorphism groups of a rapidity curve (\req(Wk'k)) for all $N$. \par \noindent
{\bf Definition}: The universal CP (chiral Potts) group ${\goth G}$ is  the group generated by ${\bf u}_1 , {\bf u}_2, {\bf U}$ with elements ${\bf U}^m {\bf u}_2^{n_2} {\bf u}_1^{n_1}$ for $m, n_1, n_2 \in \ZZ$, whose generators satisfy the relations:
\bea(ll)
(i) & {\bf u}_1{\bf u}_2 {\bf u}_1^{-1} = {\bf U}^2 {\bf u}_2^{-3} , \\
(ii) & {\bf u}_i {\bf U} {\bf u}_i^{-1} = {\bf u}_i^{-1} {\bf U} {\bf u}_i = {\bf U}^{-1} {\bf u}_i^2  ~ ~ (i=1,2).
\elea(CPG)
$\Box$ \par \vspace{.1in} \noindent
The conditions (\req(CPG)) are equivalent to the following relations:
\bea(llll)
(\req(CPG)) (ii) & \Leftrightarrow &   {\bf u}_i {\bf U} {\bf u}_i^{-1} = {\bf U}^{-1} {\bf u}_i^2  , ~ ~ ~  {\bf u}_i^2 {\bf U} = {\bf U} {\bf u}_i^2 , \\
& \Leftrightarrow &
{\bf u}_i {\bf U}^k {\bf u}_i^{-1} = {\bf U}^{-k} {\bf u}_i^{2k} ~ ~ ( k \in \ZZ ), &
\Leftrightarrow  {\bf u}_i^{-1} {\bf U}^k {\bf u}_i = {\bf U}^{-k} {\bf u}_i^{2k} ~ ~ ( k \in \ZZ ), \\
(\req(CPG)) (i) & \Leftrightarrow & {\bf u}_1 {\bf u}_2 = {\bf u}_2 {\bf U}^{-2}  {\bf u}_1 , 
\elea(CPG2)
where $i=1, 2$. By (\req(CPG2)), one obtains
\bea(llll)
{\bf u}_1 {\bf u}_2^2 {\bf u}_1^{-1} = {\bf u}_2^{-2}, & {\bf u}_1^2 {\bf u}_2^2  = {\bf u}_1^2{\bf u}_2^2 , \\
{\bf u}_2 {\bf u}_1 {\bf u}_2^{-1} = {\bf U}^{-2}{\bf u}_2^4 {\bf u}_1 , & {\bf u}_2 {\bf u}_1^2 {\bf u}_2^{-1}= {\bf u}_1^{-2} .
\elea(CPG3)
Hence ${\goth G}$ is a solvable group in which ${\bf u}_1^2, {\bf u}_2^2, {\bf U}$ generate a abelian normal subgroup of ${\goth G}$ with the quotient group $\ZZ_2^2$:
\bea(lll)
{\goth G}_1 = \langle {\bf u}_1^2, {\bf u}_2^2, {\bf U} \rangle ( \simeq \ZZ^3)  \lhd {\goth G} & \longrightarrow & \overline{\goth G}: = {\goth G}/{\goth G}_1 = \langle \overline{\bf u}_1, \overline{\bf u}_2 \rangle ( \simeq \ZZ_2^2 ).
\elea(G1)
In general, we define  \par \noindent
{\bf Definition}: For $N \geq 2$, the ($N$-state) CP group $G_N$ is the quotient group of ${\goth G}$ by the abelian (normal) subgroup ${\goth G}_N$   generated by ${\bf u}_1^{2N}, {\bf u}_2^{2N}, {\bf U}^N$:
\bea(lll)
{\goth G}_N = \langle {\bf u}_1^{2N}, {\bf u}_2^{2N}, {\bf U}^N \rangle ( \simeq \ZZ^3)  \lhd {\goth G} & \longrightarrow &  G_N:= {\goth G}/ {\goth G}_N = \langle {\sf u}_1, {\sf u}_2, U \rangle ,
\elea(CPNg)
where ${\sf u}_i (= {\sf u}_{i N}), U (= U_N)$ are the classes of ${\bf u}_i, {\bf U}$ in $G_N$. 
$\Box$ \par \vspace{.1in} \noindent
Since ${\goth G}_N \supset {\goth G}_1$, by the projections in (\req(G1)) and (\req(CPNg)), $G_N$ is a solvable group of order $4N^3$:
\bea(llll)
G_{N, 1}= \langle u_1^2, u_2^2, U \rangle ( \simeq \ZZ_N^3)  \lhd & G_N  \longrightarrow G_N/G_{N, 1} = \langle u_1, u_2 \rangle ( \simeq \ZZ_2^2 ) , & |G_N| = 4N^3.
\elea(slovN)
where $G_{N, 1}:= {\goth G}_1/{\goth G}_N$ and $G_N/G_{N, 1} = \overline{\goth G}$.
Indeed, the group $G_N$ is characterized as the group  with three generators ${\sf u}_1 , {\sf u}_2, U$ satisfying the relations in (\req(CPG)) and the finite-order condition:
\bea(lll)
{\sf u}_1^{2N} = {\sf u}_2^{2N} = U^N = 1 .
\elea(CPN)
Furthermore, the relation ${\goth G}_{N'} \supseteq {\goth G}_N$ is equivalent to 
$$
\begin{array}{llll}
G_{N'}  \longrightarrow  G_N ,  &({\sf u}_1, {\sf u}_2, U)_{N'}  \mapsto ({\sf u}_1, {\sf u}_2, U)_N & {\rm iff} ~ N| N'. \end{array}
$$
We now describe the center of ${\goth G}$ and $G_N$:
\begin{lem}\label{lem:CentG} ${\rm Cent} ({\goth G}) = {\rm id}$, and  for $N \geq 2$, 
\bea(l)
{\rm Cent}(G_N) = \left\{\begin{array}{ll}  {\rm id},  & N : {\rm odd} ; \\
\langle {\sf u}_1^N , {\sf u}_2^N  \rangle  ~ (\simeq \ZZ_2^2)   & N : {\rm even}. \end{array} \right.
\elea(CGN)
\end{lem}
{\it Proof.} First, we consider the case ${\goth G}$ and write an element ${\bf v}$ in the center ${\rm Cent} ({\goth G})$ by ${\bf v} = {\bf U}^m {\bf u}_2^{n_2} {\bf u}_1^{n_1}$ and $n_i \equiv \epsilon_i \pmod{2}$ with $\epsilon_i=0, 1$ for $i=1, 2$. By (\req(CPG2)), ${\bf u}_1 {\bf U}^m {\bf u}_1^{-1}= {\bf U}^{-m} {\bf u}_1^{2m} $ , ${\bf u}_1 {\bf u}_2^{n_2} {\bf u}_1^{-1}= ({\bf U}^{-2} {\bf u}_2^{-3})^{\epsilon_2} {\bf u}_2^{-n_2+ \epsilon_2}$, and ${\bf u}_1^2 ({\bf U}^{-2} {\bf u}_2^{-3})^{\epsilon_2} = ({\bf U}^{-2} {\bf u}_2^{-3})^{\epsilon_2} {\bf u}_1^{2 (-1)^{\epsilon_2}}$. 
Hence 
$$
{\bf v}= {\bf u}_1 {\bf v}{\bf u}_1^{-1} = {\bf U}^{-m-2\epsilon_2} {\bf u}_2^{-n_2- 2 \epsilon_2} {\bf u}_1^{n_1 + 2m (-1)^{\epsilon_2} }, 
$$
which yields ${\bf U}^{-m-2\epsilon_2}={\bf U}^m , {\bf u}_2^{-n_2- 2 \epsilon_2}={\bf u}_2^{n_2}, {\bf u}_1^{2m }=1$, equivalently 
\bea(ll)
m=0, & {\bf U}^{2\epsilon_2}= {\bf u}_2^{2n_2+ 2 \epsilon_2} = 1. 
\elea(Cent1)
Since ${\bf U}{\bf u}_i{\bf U}^{-1}= {\bf U}^2 {\bf u}_i^{-1}$, we have 
$$
\begin{array}{lll}
{\bf u}_2^{n_2} {\bf u}_1^{n_1}= {\bf v}= {\bf U} {\bf v}{\bf U}^{-1} &= {\bf u}_2^{n_2- \epsilon_2} ( {\bf U}^2 {\bf u}_2^{-1})^{\epsilon_2} {\bf u}_1^{n_1-\epsilon_1}( {\bf U}^2 {\bf u}_i^{-1})^{\epsilon_1} & 
= U^{2 \epsilon_2+ 2 \epsilon_1(-1)^{\epsilon_2}} {\bf u}_2^{n_2- 2\epsilon_2+4 \epsilon_1\epsilon_2}   {\bf u}_1^{n_1- 2\epsilon_1}, 
\end{array}
$$
equivalently 
\bea(ll)
\epsilon_1 = \epsilon_2=0, &i.e. ~ n_1 \equiv n_2 \equiv 0  \pmod{2}.
\elea(Cent2)
Then ${\bf u}_2 {\bf v} {\bf u}_2^{-1} = {\bf u}_2^{n_2} {\bf u}_1^{-n_1}$, i.e. 
\be
{\bf u}_1^{2n_1} = 1 .
\ele(Cent3)
Then by (\req(Cent1)), (\req(Cent2)) and (\req(Cent3)), ${\bf v}=1$. In the case $G_N ~ (N \geq 2)$, (\req(Cent1)), (\req(Cent2)) and (\req(Cent3))  again hold for an element 
${\sf v} \in {\rm Cent}(G_N)$, i.e. ${\sf v} =  {\sf u}_2^{n_2} {\sf u}_1^{n_1}$ with $n_i \equiv 0  \pmod{2}$ and 
${\bf u}_i^{2n_i} = 1$ for $i=1, 2$. Hence follows the relation (\req(CGN)). 
$\Box$ \par \vspace{.2in} \noindent
For an element ${\bf v} \in {\goth G}$, the conjugation of ${\bf v}$ on ${\goth G}$ will be denoted by
\bea(ll)
{\sf C}_{\bf v} : {\goth G} \longrightarrow {\goth G} & {\bf g} \mapsto {\bf v} {\bf g} {\bf v}^{-1},
\elea(conjG)
which preserves the normal subgroup ${\goth G}_1$ in (\req(G1)), and induces the identity of the quotient group $\overline{\goth G}$. Furthermore, ${\goth G}_N$  in (\req(CPNg)) is preserved by ${\sf C}_{\bf v}$, which induces the conjugation of ${\sf v} \in G_N$ on $G_N$: 
\bea(ll)
C_{\sf v} : G_N \longrightarrow G_N & g \mapsto {\sf v} g {\sf v}^{-1}.
\elea(conjGN)
By Lemma \ref{lem:CentG}, ${\goth G}$ can regarded as a subgroup of the automorphism group ${\rm Aut}({\goth G})$ of ${\goth G}$ by the conjugation (\req(conjG)):
\bea(ll)  
{\sf C}: {\goth G} \hookrightarrow {\rm Aut}({\goth G}) , & {\bf v} \mapsto {\sf C}_{\bf v},
\elea(CG)
and (\req(conjGN)) defines an embedding of $G_N/{\rm Cent}(G_N)$ in ${\rm Aut}(G_N)$:
\bea(lll)  
{\rm Cent}(G_N) \hookrightarrow & G_N  \stackrel{C}{\rightarrow}& {\rm Aut}(G_N). 
\elea(CGN) 
Using (\req(CPG2)), one finds the relations $({\bf U}^{-1} {\bf u}^{-3}_1 {\bf U} ){\bf u}_2 = {\bf u}_2^{-3} {\bf u}_1^3$ and ${\bf u}_1 = {\bf u}_2 {\bf U}^{-2} ( {\bf u}_1 {\bf u}_2^{-1})$. By which, the following correspondences of generators give rise to two automorphisms of ${\goth G}$:
\bea(rll)
{\bf S} : {\goth G} \longrightarrow {\goth G}, &({\bf u}_1, {\bf u}_2, {\bf U} ) \mapsto ({\bf u}_2 , {\bf U}^{-2}{\bf u}_1^3, {\bf U})=
({\bf u}_2 , {\bf U}^{-1}{\bf u}_1 {\bf U}, {\bf U}) , \\
( \Leftrightarrow {\bf S}^{-1}:& 
({\bf u}_1, {\bf u}_2, {\bf U}) \mapsto ( {\bf U}^2{\bf u}_2^{-1}, {\bf u}_1 , {\bf U}) )= ( {\bf U}{\bf u}_2{\bf U}^{-1}, {\bf u}_1 , {\bf U}) );  \\
{\bf T} : {\goth G} \longrightarrow {\goth G}, &({\bf u}_1, {\bf u}_2, {\bf U} ) \mapsto ( {\bf U}^2 {\bf u}_2^{-1}{\bf u}_1, {\bf u}_2 , {\bf U})= ({\bf u}_1 {\bf u}_2^{-1} , {\bf u}_2, {\bf U}),  \\
( \Leftrightarrow {\bf T}^{-1}:&({\bf u}_1, {\bf u}_2, {\bf U}) \mapsto ( {\bf U}^2 {\bf u}_2^{-3}{\bf u}_1, {\bf u}_2, {\bf U} )=( {\bf u}_1{\bf u}_2, {\bf u}_2, {\bf U} ) )
\elea(ST)
which satisfy the relations:
\bea(lll)
{\bf S}^2 = ({\bf S} {\bf T})^3 = {\sf C}_{{\bf U}^{-1}} , & {\bf T}^4 = {\sf C}_{{\bf u}_2^2}. 
\elea(PowST)
The first equality in (\req(PowST)) is equivalent to ${\bf S}{\bf T}{\bf S}^{-1}= ({\bf T}{\bf S}{\bf T})^{-1}$, which implies 
$$
{\bf S}^2{\bf T}{\bf S}^{-2}= ({\bf S}{\bf T}^{-1}{\bf S}^{-1}){\bf T}^{-1}{\bf S}^{-1}= {\bf T}, ~ ~ {\rm i.e.} ~ ~ {\bf S}^2 {\bf T} = {\bf T}  {\bf S}^2.
$$
Indeed by ${\bf S}^2 = {\sf C}_{{\bf U}^{-1}}$, ${\bf S}^2 {\bf T} = {\bf T}  {\bf S}^2$ ( or  
${\bf T}^{-1}{\bf S}^2 {\bf T} = {\bf S}^2$) is equivalent to ${\sf C}_{{\bf T}^{-1}({\bf U}^{-1})}= {\sf C}_{{\bf U}^{-1}}$ (or ${\bf T} ({\bf U}) = {\bf U} $). Since the automorphisms in (\req(ST)) preserve ${\goth G}_1, {\goth G}_N$ in (\req(G1)) and (\req(CPNg)) by
\bea(ll)
{\bf S} : ({\bf u}_1^2, {\bf u}_2^2, {\bf U} ) \mapsto ({\bf u}_2^2 , {\bf u}^2_1, {\bf U}), \\
{\bf T} : ({\bf u}^2_1, {\bf u}^2_2, {\bf U} ) \mapsto ({\bf U}^2 {\bf u}_2^{-2}{\bf u}_1^{-2} , {\bf u}^2_2, {\bf U}),
\elea(STG1)
${\bf S}, {\bf T}$ induce the automorphisms of $\overline{\goth G}$ and $G_N$:
\bea(lll)
\overline{\bf S}, \overline{\bf T}  : \overline{\goth G} \longrightarrow \overline{\goth G}, &\overline{\bf S} (\overline {\bf u}_1, \overline {\bf u}_2) = (\overline {\bf u}_2 , \overline {\bf u}_1), & \overline{\bf T} (\overline {\bf u}_1, \overline {\bf u}_2)=  (\overline {\bf u}_1 \overline {\bf u}_2 , \overline {\bf u}_2) ; \\
{\sf S} , {\sf T}  : G_N \longrightarrow G_N , & {\sf S} ({\sf u}_1, {\sf u}_2, U ) =
({\sf u}_2 , U^{-1}{\sf u}_1 U, U) , & {\sf T} ({\sf u}_1, {\sf u}_2, U)= ({\sf u}_1 {\sf u}_2^{-1} , {\sf u}_2, U).
\elea(STquot)
The automorphisms in (\req(ST)), (\req(STquot)) then generate the automorphism subgroups of ${\goth G}, G_N$ , denoted respectively by  
$$
\langle {\bf S}, {\bf T} \rangle \subseteq {\rm Aut}({\goth G}), ~ \langle {\sf S}, {\sf T} \rangle \subseteq {\rm Aut}(G_N), ~  \langle \overline{\bf S}, \overline{\bf T} \rangle \subseteq {\rm Aut}(\overline{\goth G}), ~ ~ (N \geq 2),
$$
compatible with the projections of ${\goth G}$ to $G_N, \overline{\goth G}$ in (\req(CPNg)), (\req(slovN)) respectively:
\bea(ccccc)
\langle {\bf S}, {\bf T} \rangle & \longrightarrow & \langle {\sf S}, {\sf T} \rangle & \longrightarrow &\langle \overline{\bf S}, \overline{\bf T} \rangle 
\elea(STN1)
With the identification $\overline{\goth G}= \ZZ^2$ in (\req(G1)), one finds $\langle \overline{\bf S}, \overline{\bf T} \rangle = SL_2(\ZZ_2)$:
\bea(ll)
\overline{\bf S} ~ (= \overline{S}) = \left( \begin{array}{cc}
0&1\\ 
1&0
\end{array}
\right), & \overline{\bf T} ~ (= \overline{T}) = \left( \begin{array}{cc}
1& 0 \\ 
1 &1
\end{array} \right) \in SL_2(\ZZ_2) .
\elea(ST-) 
\begin{prop}\label{prop:Ker}
The correspondence of $\langle {\bf S}, {\bf T} \rangle$ to $\langle {\sf S}, {\sf T} \rangle$ induces the
canonical isomorphisms:
\bea(llll)
\langle {\bf S}, {\bf T} \rangle/ (\langle {\bf S}, {\bf T} \rangle \cap {\sf C}({\goth G})) & \simeq & \langle {\sf S}, {\sf T} \rangle/ (\langle {\sf S}, {\sf T} \rangle \cap C(G_N))  & \simeq PSL_2(\ZZ_4)
\elea(PSLZ4)
with $(\widetilde{\bf S}, \widetilde{\bf T}) \leftrightarrow (\widetilde{{\sf S}}, \widetilde{\sf T}) \leftrightarrow ({\rm S}, {\rm T}^*) $ satisfying the relations: 
\be
\widetilde{\bf S}^2 = (\widetilde{\bf S} \widetilde{\bf T})^3 = \widetilde{\bf T}^4 = 1 ,
\ele(STrel)
where ${\rm S}, {\rm T}^*$ are the standard generators of  $PSL_2(\ZZ_4)$:
\bea(ll)
{\rm S} = \left( \begin{array}{cc}
0& -1\\ 
1&0
\end{array}
\right), & {\rm T}^* = \left( \begin{array}{cc}
1& 0 \\ 
-1 &1
\end{array}
\right).
\elea(STSL4)
\end{prop}
{\it Proof.} Note that ${\bf u}_i{\bf U}{\bf u}^{-1}_i ={\bf U}^{-1} {\bf u}_i^2 \neq {\bf U}$ and ${\bf u}_i{\bf u}_j{\bf U}{\bf u}^{-1}_j {\bf u}_i^{-1} ={\bf U}{\bf u}_2^2 {\bf u}_1^2 \neq {\bf U}$, in ${\goth G}$ (or $G_N$ by $N \geq 2$). Since both ${\bf S}, {\bf T}$ fix the element ${\bf U}$, $\langle {\bf S}, {\bf T} \rangle \cap {\sf C}({\goth G}) \subseteq {\sf C}({\goth G}_1)$ where ${\goth G}_1$ is the abelian subgroup in (\req(G1)). Indeed, by (\req(ST)), (\req(PowST)) and $({\bf S}^{-1} {\bf T} {\bf S})^4 = {\sf C}_{{\bf u}_1^2}$, one finds  
\bea(llll)
{\sf C}({\goth G}) \cap \langle {\bf S}, {\bf T} \rangle &= \langle \langle {\bf S}^2, ({\bf S} {\bf T})^3 , {\bf T}^4 \rangle \rangle_{\rm N} = \langle {\bf S}^2,  {\bf T}^4, ({\bf S}^{-1} {\bf T} {\bf S})^4 \rangle = {\sf C}({\goth G}_1)  & \lhd   \langle {\bf S}, {\bf T} \rangle ; \\
C(G_N) \cap \langle {\sf S}, {\sf T} \rangle  &= \langle \langle {\sf S}^2, ({\sf S}{\sf T})^3 , {\sf T}^4 \rangle \rangle_{\rm N} = \langle {\sf S}^2, {\sf T}^4, ({\sf S}^{-1} {\sf T} {\sf S})^4 \rangle = C(G_{N, 1}) &\lhd  \langle {\sf S}, {\sf T} \rangle 
\elea(STCG1)
where $\langle \langle * \cdots * \rangle \rangle_{\rm N}$ is the normal subgroup defined in (\req(BNH)). Hence the relation (\req(STrel)) holds. By a well-known characterization of $PSL_2(\ZZ_4)$, $\langle \widetilde{\bf S}, \widetilde{\bf T} \rangle$ is a quotient group of $PSL_2(\ZZ_4)$ by assigning ${\rm S}, {\rm T}^*$ 
in $PSL_2(\ZZ_4)$ to  $\widetilde{\bf S}, \widetilde{\bf T}$ respectively. Together with the induced homomorphisms of  quotient groups in (\req(STN1)), we obtain the following group epimorphisms of quotient groups:
 $$
\begin{array}{lllll}
PSL_2(\ZZ_4) \longrightarrow \langle \widetilde{\bf S}, \widetilde{\bf T} \rangle & \longrightarrow & \langle \widetilde{\sf S}, \widetilde{\sf T}  \rangle & \longrightarrow &\langle \overline{\bf S}, \overline{\bf T} \rangle = SL_2(\ZZ_2).
\end{array}
$$
Since the kernel of the homomorphism from $PSL_2(\ZZ_4)$ to $SL_2(\ZZ_2)$ is equal to $\langle ({\rm T}^{* 2}, {\rm S} {\rm T}^{* 2} S_4^{-1} \rangle (\simeq \ZZ_2^2)$, which contains no proper non-trivial normal subgroup of $PSL_2(\ZZ_4)$,  the projection from $PSL_2(\ZZ_4)$ to $\langle \widetilde{\sf S}, \widetilde{\sf T}  \rangle$ must have the trivial kernel by $\widetilde{\sf T} ^2 \neq 1$, then follow the isomorphisms in (\req(PSLZ4)).
$\Box$ \par \vspace{.1in} \noindent
 {\bf Remark.} The $PSL_2(\ZZ_4)$-structure of $\langle \widetilde{\sf S}, \widetilde{\sf T}  \rangle$ in (\req(PSLZ4)) for $N=2$ can be explicitly derived by the group structure of $G_2$, where the conditions, (\req(CPG)) and (\req(CPN)), for the $G_2$-generators ${\sf u}_1, {\sf u}_2, U$  are equivalent to 
\bea(lll)
{\sf u}_1^4 = {\sf u}_2^4 = U^2=1, & {\sf u}_1{\sf u}_2= {\sf u}_2{\sf u}_1, & U{\sf u}_i U^{-1} = {\sf u}_i^{-1} ~ ~ (i=1,2). 
\elea(G2)
Hence $\langle {\sf u}_1, {\sf u}_2 \rangle $ is a normal abelian subgroup of $G_2$. By (\req(STquot)) and (\req(G2)), both $\langle {\sf S}, {\sf T} \rangle$ and $C({\goth G}_1)$ fix $U$, and leave $\langle {\sf u}_1, {\sf u}_2 \rangle (\simeq \ZZ_4^2)$ invariant, with the generators represented by the following elements in $SL_2(\ZZ_4)$:
$$
({\sf S}, {\sf T}, C_U, C_{{\sf u}^2_1}, C_{{\sf u}^2_2}) \leftrightarrow  ( {\rm S},   {\rm T}^* ,   -{\rm id.} , {\rm id.}, ~ {\rm id.} ). 
$$
By (\req(PowST)), $\langle {\sf S}, {\sf T} \rangle \cap {\sf C}({\goth G}_1)$ is generated by the $-{\rm id.}$ in $SL_2(\ZZ_4)$, hence $\langle \widetilde{\sf S}, \widetilde{\sf T}  \rangle \simeq PSL_2(\ZZ_4)$ in (\req(PSLZ4)). 
$\Box$ \par \vspace{.1in}  
Consider the semi-product of ${\goth G}$ and $\langle {\bf S}, {\bf T} \rangle$, ${\goth G} \ast \langle {\bf S}, {\bf T} \rangle$, with the group-multiplication 
$$
\begin{array}{lll}
({\bf v} \ast {\bf M}) \cdot({\bf v}' \ast {\bf M}') := {\bf v}{\bf M}({\bf v}') \ast {\bf M} {\bf M}' & {\rm for} ~  {\bf M}, {\bf M}' \in \langle {\bf S}, {\bf T} \rangle, & {\bf v}, {\bf v}' \in {\goth G}, 
\end{array}
$$
and define  \par \noindent
{\bf Definition}: The universal CP modular group 
\bea(ll)
\widetilde{\goth G} : = \bigg({\goth G} \ast \langle {\bf S}, {\bf T} \rangle \bigg)/\langle \langle {\bf U}* {\bf S}^2, {\bf U}*({\bf S}{\bf T})^3 , {\bf u}_2^{-2}*{\bf T}^4  \rangle \rangle_{\rm N}. 
\elea(uCPm)
$\Box$ \par \vspace{.1in} \noindent
\begin{lem}\label{lem:NomCPM}
The normal subgroup in (\req(uCPm)) is isomorphic to ${\goth G}_1$ in (\req(G1)) by: 
\bea(ll)
\langle \langle {\bf U}* {\bf S}^2, {\bf U}*({\bf S}{\bf T})^3 , {\bf u}_2^{-2}*{\bf T}^4  \rangle \rangle_{\rm N} = \{ {\sf u}^{-1}* {\sf C}_{\sf u} \in {\goth G} \ast \langle {\bf S}, {\bf T} \rangle | ~ {\sf u} \in {\goth G}_1 \}.
\elea(NuCPM) 
\end{lem}
{\it Proof.} In ${\goth G} \ast \langle {\bf S}, {\bf T} \rangle$, one finds
$$
\begin{array}{ll}
({\sf u}^{-1}* {\sf C}_{\sf u})\ast ({\sf u}'^{-1}* {\sf C}_{\sf u'} ) = ({\sf u u'})^{-1}*{\sf C}_{\sf u u'},   \\
{\sf v} \ast ({\sf u}^{-1}* {\sf C}_{\sf u}) \ast {\sf v}^{-1} = {\sf u}^{-1} \ast {\sf C}_{\sf u}, &
{\bf M} \ast ({\sf u}^{-1}* {\sf C}_{\sf u}) \ast {\bf M}^{-1} = {\bf M}({\sf u})^{-1}* {\sf C}_{{\bf M} ({\sf u})},
\end{array}
$$
for ${\sf u}, {\sf u}' \in {\goth G}_1 $, ${\sf v} \in {\goth G}, {\bf M} \in \langle {\bf S}, {\bf T} \rangle$. Hence the right hand side of (\req(NuCPM)) is an abelian normal subgroup of ${\goth G} \ast \langle {\bf S}, {\bf T} \rangle$. 
Then (\req(NuCPM)) follows from the following equalities:   
$$
\begin{array}{lll}
{\bf U}* {\bf S}^2= {\bf U}*({\bf S}{\bf T})^3 = {\bf U}*{\sf C}_{{\sf U}^{-1}} ,& {\bf u}_2^{-2}*{\bf T}^4={\bf u}_2^{-2}*{\sf C}_{{\bf u}_2}, & {\bf u}_1^{-2}*({\bf S}^{-1}{\bf T}{\bf S})^4={\bf u}_1^{-2}*{\sf C}_{{\bf u}_1}.
\end{array}
$$
$\Box$ \par \vspace{.1in} \noindent
By (\req(STCG1)), (\req(NuCPM)),  with the identification of ${\goth G}$ and the inner automorphism group of ${\goth G}$ via the conjugation (\req(CG)), $\widetilde{\goth G}$ in (\req(uCPm)) can be regarded as an automorphism group of ${\goth G}$:
\be
\widetilde{\goth G} \simeq \langle {\sf C} ({\goth G}), {\bf S}, {\bf T} \rangle \subset {\rm Aut}({\goth G}). 
\ele(UGtAut)
By Lemma \ref{lem:CentG} and (\req(NuCPM)), ${\goth G}$ and $\langle {\bf S}, {\bf T} \rangle$ are embedded as subgroups of $\widetilde{\goth G}$. For convenience, we shall write the class of ${\sf v} \in {\goth G}$ again by ${\sf v}$, and the class of $1 \ast {\bf S}, 1 \ast {\bf T}$ or $1 \ast {\bf M} \in  \langle {\bf S}, {\bf T} \rangle$ in $\widetilde{\goth G}$ by $S, T$ or $M$ respectively. Note that  $M{\bf v}M^{-1} = {\bf M}({\bf v})$ in $\widetilde{\goth G}$, so ${\bf M}$ is identified with ${\sf C}_M := \widetilde{\sf C}_{M | {\goth G}}$, the restriction of $\widetilde{\sf C}_M$ on ${\goth G}$, where $\widetilde{\sf C}_{\bf g}:= $ the ${\bf g}$-conjugation of $\widetilde{\goth G}$ for ${\bf g} \in \widetilde{\goth G}$. By (\req(PSLZ4)), (\req(NuCPM)) and  (\req(UGtAut)), one obtains
\bea(llll)
{\goth G} \lhd \widetilde{\goth G}= {\goth G}  \langle S, T \rangle , & S^2 =  (ST)^3 = {\bf U}^{-1}, &T^4= {\bf u}_2^2 ;& \widetilde{\goth G} / {\goth G}  \simeq PSL_2(\ZZ_4); \\
{\goth G}  \cap \langle S, T \rangle = {\goth G}_1 &=  \langle \langle S^2, (ST)^3 , T^4 \rangle \rangle_{\rm N} &=
\langle S^2,  T^4, (S^{-1} T S)^4 \rangle &\lhd   \langle S, T \rangle.
\elea(uGtil)
Indeed, $\widetilde{\goth G}$ is defined the first three properties in (\req(uGtil)) as follows:
\begin{prop}\label{prop:utilG}
$\widetilde{\goth G}$ is charactered as the group generated by ${\goth G}$ and $S, T$ satisfying the relations: 
\bea(llll)
S {\bf v} S^{-1} = {\bf S}({\bf v}), & T {\bf v} T^{-1} = {\bf T}({\bf v}),
 & S^2 =  (ST)^3 = {\bf U}^{-1}, & T ^4= {\bf u}_2^2 ,  
\elea(uCPrel)
where ${\bf v} \in {\goth G}$, and ${\bf S}, {\bf T}$ are defined in (\req(ST)).
\end{prop}
{\it Proof.} Let $\langle {\goth G}, S, T \rangle$ be the group generated by ${\goth G}$ and $S, T$ defined by the relation (\req(uCPrel)). Then there is a group epimorphism, $\wp: \langle {\goth G}, S, T \rangle \longrightarrow \widetilde{\goth G}$ which is the identity on the normal subgroup ${\goth G}$. By the characterization of $PSL_2(\ZZ_4)$, $\wp$ induces an isomorphism between $\langle {\goth G}, S, T \rangle/ {\goth G}$ and $\widetilde{\goth G} / {\goth G}$, hence ${\rm Ker}(\wp) \subset {\goth G} \subset \langle {\goth G}, S, T \rangle$. Then $\wp$ defines the isomorphism between $\langle {\goth G}, S, T \rangle$ and $\widetilde{\goth G}$.
$\Box$ \par \vspace{.1in} 
By (\req(STG1)), the abelian normal subgroups ${\goth G}_1, {\goth G}_N$ of ${\goth G}$ in (\req(G1)) (\req(CPNg)) are also normal in $\widetilde{\goth G}$. \par  \noindent
{\bf Definition}: For $N \geq 2$, the modular ($N$-state) CP group is the quotient group of $\widetilde{\goth G}$ by 
${\goth G}_N$:
\bea(lll)
{\goth G}_N  \lhd \widetilde{\goth G} & \longrightarrow &  \widetilde{G}_N:= \widetilde{\goth G}/ {\goth G}_N = \langle G_N, S, T \rangle ,
\elea(CPNmg)
where $S (= S_N), T (=T_N)$ are the class of the $\widetilde{\goth G}$-elements $S, T$ in $\widetilde{G}_N$. $\Box$ \par \vspace{.1in}  \noindent
By (\req(slovN)), (\req(uGtil)) and Proposition \ref{prop:utilG}, we obtain the following results:
\begin{prop}\label{prop:CPMm} (i) $\widetilde{G}_N$ is charactered as the group generated by $G_N$ and $S, T$ satisfying the relations: 
\bea(llll)
S {\sf v} S^{-1} = {\sf S}({\sf v}), & T {\sf v} T^{-1} = {\sf T}({\sf v}),
 & S^2 =  (ST)^3 = U^{-1}, & T ^4= {\sf u}_2^2 ,  
\elea(tCPNrel)
where ${\sf v} \in G_N$, and ${\sf S}, {\sf T}$ are defined in (\req(STquot)).

(ii) $G_N$ and $G_{N, 1} (= {\goth G}_1/{\goth G}_N)$ are normal subgroups of $\widetilde{G}_N$ with the following relations:
\bea(llll)
G_N \lhd \widetilde{G}_N= G_N  \langle S, T \rangle &\longrightarrow  \widetilde{G}_N / G_N  \simeq PSL_2(\ZZ_4), &|\widetilde{G}_N| = 96 N^3 ; \\
G_N  \cap \langle S, T \rangle = G_{N, 1} = & \langle \langle S^2, (ST)^3 , T^4 \rangle \rangle_{\rm N} =
\langle S^2,  T^4, (S^{-1} T S)^4 \rangle &\lhd   \langle S, T \rangle.
\elea(tGNrel)
\end{prop}
$\Box$ \par \vspace{.2in} \noindent
For an element ${\sf V}\in \widetilde{G}_N$, we denote  the ${\sf V}$-conjugation of $\widetilde{G}_N$ by $\widetilde{C}_{{\sf V}}$, and its  restriction on  $G_N$ by 
$$
C_{\sf V} := \widetilde{C}_{{\sf V} |G_N}: G_N \longrightarrow G_N , ~ ~ g \mapsto {\sf V} g {\sf V}^{-1}. 
$$
which is the same as $C_{\sf v}$ in (\req(conjGN)) for ${\sf v} \in G_N$. The first two relations in (\req(tCPNrel)) mean 
\bea(ll)
C_S = {\sf S}, & C_T = {\sf T}, 
\elea(CgST) 
by which, $M \in \langle S, T \rangle \subset \widetilde{G}_N$ corresponds to ${\sf M} \in \langle {\sf S}, {\sf T} \rangle$ with the relation $C_M = {\sf M}$. By (\req(PSLZ4)), the conjugation morphism
\bea(ll)
C: \widetilde{G}_N \longrightarrow C(G_N)\langle {\sf S}, {\sf T} \rangle \subset {\rm Aut}(G_N) , ~ ~ {\sf V} \mapsto C_{\sf V} ,
\elea(conjtGN)
gives rise to the isomorphism
\bea(ll)
\widetilde{G}_N/G_N \simeq \langle {\sf S}, {\sf T} \rangle/ (\langle {\sf S}, {\sf T} \rangle \cap C(G_N))   \simeq PSL_2(\ZZ_4). 
\elea(Gt/GN)
Since ${\rm Cent}(PSL_2(\ZZ_4)) = 1$, the kernel of $C$ in (\req(conjtGN)) is contained in $G_N$, and by (\req(CGN)),
\bea(ll)
1 = {\rm Cent}(\widetilde{G}_N) \subset {\rm Ker}(C) = {\rm Cent}(G_N).
\elea(CentGN)
Indeed for even $N$, by (\req(STG1)) ${\sf S}({\sf u}_1^N, {\sf u}_2^N)= ({\sf u}_2^N, {\sf u}_1^N)$ and 
${\sf T}({\sf u}_1^N, {\sf u}_2^N)= ({\sf u}_1^N {\sf u}_2^N, {\sf u}_2^N)$, hence ${\rm Cent}(\widetilde{G}_N) = \{ {\sf v} \in {\rm Cent}(G_N) | S {\sf v} S^{-1} = T {\sf v} T^{-1} = {\sf v} \}= 1$. Note that by composing with the morphism in (\req(STN1)), the projection in (\req(Gt/GN)) gives rise to the isomorphism
\bea(ll)
\widetilde{G}_N/\langle G_N, T^2, ST^{-2}S^{-1} \rangle \stackrel{\sim}{\longrightarrow}  SL_2(\ZZ_2), & ({\sf v}, S, T ) \mapsto (1, \overline{S}, \overline{T}), 
\elea(GNSL22)
where ${\sf v} \in G_N$ and $\overline{S}, \overline{T}$ are defined in (\req(ST-)).

\subsection{CP group and $\ZZ_2 \times D_N$ \label{ssec:Dih}}
The spectral parameter of $\tau^{(2)}$-model in CPM lies in a hyperelliptic curves with $(\ZZ_2 \times D_N)$-symmetry, where  $D_N$ is the dihedral group. In this subsection, we study the $(\ZZ_2 \times D_N)$-structures reduced from the CP group $G_N$ from the group-theory point of view, and examine their relationship under the action of $\widetilde{G}_N$.
First, we consider the subgroup of ${\goth G}_1$ in (\req(G1)), which is normal in the universal CP group $\widetilde{\goth G}$:
\bea(lll)
{\goth G}'_1 = \langle {\bf V}_0, {\bf V}_1,  {\bf V}_2 \rangle & ( \simeq \ZZ^3) \lhd  \widetilde{\goth G}, & {\goth G}'_1 \subset {\goth G}_1, \\
{\bf V}_0:={\bf U}^2 {\bf u}_2^{-2}{\bf u}_1^{-2}, &{\bf V}_1:= {\bf u}_1^2, & {\bf V}_2:= {\bf u}_2^2,
\elea(G12)
where the $\widetilde{\goth G}$-conjugation on the generators of ${\goth G}'_1$ is given by ${\sf C}_{{\bf U}|{\goth G}'_1}= {\rm id}$, and  
\bea(cll)
{\sf C}_{{\bf u}_1}, {\sf C}_{{\bf u}_2}&: ({\bf V}_0, {\bf V}_1,  {\bf V}_2) &\mapsto ({\bf V}_0^{-1},  {\bf V}_1, {\bf V}_2^{-1} ), ({\bf V}_0^{-1}, {\bf V}_1^{-1},  {\bf V}_2 ); \\
{\bf S}, {\bf T}  &: ({\bf V}_0,  {\bf V}_1,  {\bf V}_2) & \mapsto ({\bf V}_0,  {\bf V}_2 ,  {\bf V}_1 ), ( {\bf V}_1 , {\bf V}_0 ,  {\bf V}_2).  
\elea(STG12)
For $N \geq 2$, the image of ${\goth G}'_1 $ in $G_N$ under the projection in (\req(CPNmg)) will be denoted by  
\bea(llll)
G'_{N, 1} = \langle {\sf V}_0, {\sf V}_1, {\sf V}_2  \rangle &\lhd  \widetilde{G}_N, & G'_{N, 1} \subset G_{N, 1}, &
{\sf V}_0:=U^2 {\sf u}_2^{-2}{\sf u}_1^{-2}, ~ {\sf V}_1:= {\sf u}_1^2, ~ {\sf V}_2:= {\sf u}_2^2,
\elea(GN12)
\begin{lem}\label{lem:GG12} The quotient group $\widetilde{G}_N/G'_{N, 1}$ is given by 
\bea(lll)
\widetilde{G}_N/G'_{N, 1} = \left\{ \begin{array}{lll} \widetilde{\goth G}/{\goth G}'_1 \simeq \ZZ_2^2 \ast SL_2(\ZZ_4),& N: {\rm even} , \\
\widetilde{\goth G}/{\goth G}_1 \simeq \ZZ_2^2 \ast PSL_2(\ZZ_4), & N: {\rm odd}. 
\end{array} \right.
\elea(G/G12)
\end{lem}
{\it Proof.} When $N$ is even, ${\goth G}_N$ of ${\goth G} $ in (\req(CPNg)) is a subgroup  of ${\goth G}'_1$, hence $\widetilde{G}_N/G'_{N, 1} = \widetilde{\goth G}/{\goth G}'_1$. By (\req(PowST)) and Proposition \ref{prop:CPMm} (i) , $\widetilde{\goth G}/{\goth G}'_1$ is the group generated by
$[{\bf u}_1], [{\bf u}_2], [S], [T]$ with the relations 
\bea(lll)
&[{\bf u}^2_1]= [{\bf u}^2_2]=  {\rm id}, & [{\bf u}_1{\bf u}_2]=[ {\bf u}_2{\bf u}_1] , \\
&[S ({\bf u}_1, {\bf u}_2) S^{-1}] = [({\bf u}_2, {\bf u}_1)], &[T ({\bf u}_1, {\bf u}_2)T^{-1}] = [({\bf u}_1 {\bf u}_2, {\bf u}_2)]
\elea(Z2ST)
$[S^4]= [T^4] = {\rm id}$ and $[S^2] = [S T]^3$, which characterize the semi-product group $\ZZ_2^2 \ast SL_2(\ZZ_4)$. For odd $N$, ${\goth G}_N$ is not contained in ${\goth G}'_1$. Indeed in this case,  ${\goth G}_N$ and ${\goth G}'_1$ generate the subgroup ${\goth G}_1$ in (\req(G1)), hence $\widetilde{G}_N/G'_{N, 1} = 
\widetilde{\goth G}/{\goth G}_1$,  which is isomorphic to $\ZZ_2^2 \ast PSL_2(\ZZ_4)$ with the generators 
$[{\bf u}_1], [{\bf u}_2], [S], [T]$, satisfying the relation (\req(Z2ST)) and $[S^2]= [S T]^3 = [T^4] = {\rm id}$.
$\Box$ \par \vspace{.2in}
We now describe another expression for the universal CP group ${\goth G}$.
Consider the following subgroups of ${\goth G}$:
\bea(lll)
{\goth H} := \langle {\bf V}_0, {\bf V}_1 \rangle, & 
{\goth H}_l := \langle {\bf V}_0, {\bf V}_2 \rangle, & {\goth H}_r := \langle {\bf V}_1, {\bf V}_2 \rangle.
\elea(uHrl)
By (\req(STG12)), the above groups are normal in ${\goth G}$, interchanged under $\langle {\bf S}, {\bf T} \rangle$. First, we give another characterization of ${\goth G}$ through the abelian normal subgroup ${\goth H}$ in (\req(uHrl)). 
\begin{lem}\label{lem:Uij} The universal CP group ${\goth G}$ is characterized as the group with generators ${\bf V}_0, {\bf V}_1, {\bf U}$ and ${\bf i}, {\bf j}$  satisfying the relations:
\bea(llll)
{\bf V}_0 {\bf V}_1 = {\bf V}_1{\bf V}_0, & {\sf C}_{{\bf U}}({\bf V}_k) = {\bf V}_k, & {\sf C}_{{\bf j}}({\bf V}_k) = {\bf V}_k^{-1} &{\sf C}_{{\bf i}}({\bf V}_k) = {\bf V}_k^{(-1)^k}, ~ ~  (k=0, 1),  \\
 {\bf j}^2 = {\bf i}^2  = 1, &{\bf U}{\bf j} = {\bf V}_0 {\bf V}_1 {\bf j} {\bf U}, & {\bf U}{\bf i} = {\bf V}_0 {\bf i} {\bf U}^{-1}, & {\bf j} {\bf i} = {\bf V}_1 {\bf i} {\bf j}, 
\elea(bVUij)
where ${\bf V}_0, {\bf V}_1$ and ${\bf i}, {\bf j}$ are related to the generators ${\bf U}, {\bf u}_1, {\bf u}_2$ of ${\goth G}$ by the relations
\bea(llll)
{\bf V}_0= {\bf U}^2 {\bf u}_2^{-2}{\bf u}_1^{-2}, ~ {\bf V}_1= {\bf u}_1^2, & {\bf j} = {\bf U} {\bf u}_2^{-1}, ~ ~ {\bf i} = {\bf U} {\bf u}_2^{-1} {\bf u}_1, &
\Longleftrightarrow & {\bf u}_1 = {\bf j}{\bf i}, ~ {\bf u}_2 = {\bf j}{\bf U}.
\elea(u01IJ)
\end{lem}
{\it Proof.} With the expression of ${\bf V}_k, {\bf i}, {\bf j}$ in terms of  ${\bf U}, {\bf u}_1, {\bf u}_2$ in (\req(u01IJ)),  (\req(CPG)) and (\req(CPG2)) yield  the relation (\req(bVUij)). Conversely, if ${\bf V}_0, {\bf V}_1$ and ${\bf U}, {\bf i}, {\bf j}$ satisfy the relation (\req(bVUij)), using the expression of ${\bf u}_1, {\bf u}_2$ in terms of 
${\bf U}, {\bf i}, {\bf j}$, one finds all relations in (\req(CPG)) are valid. Then the last two relations in (\req(bVUij)) yield all relations in (\req(u01IJ)).
$\Box$ \par \vspace{.1in} \noindent
By Lemma \ref{lem:Uij}, ${\goth G}= \langle {\bf U},{\bf j},{\bf i} \rangle$, by which one obtains the relation between ${\goth G}$ and universal dihedral group ${\goth D}$:
\bea(lll)
{\goth H} ~ (= \langle {\bf V}_0 , {\bf V}_1 \rangle  \simeq \ZZ^2 )  \lhd & {\goth G} \longrightarrow  {\goth G}/{\goth H}  = \ZZ_2 \times {\goth D} ~ (=  \langle \Upsilon \rangle \times \langle \Theta, {\goth I} \rangle) , \\
\Upsilon^2= {\goth I}^2=1, &\Upsilon\Theta= \Theta\Upsilon, ~ ~ \Upsilon{\goth I}= {\goth I}\Upsilon, ~ ~
\Theta{\goth I} = {\goth I}\Theta^{-1}
\elea(uGDih)
where $\Theta, \Upsilon, {\goth I}$ are the classes of ${\bf U},{\bf j},{\bf i}$ in ${\goth G}/\langle {\bf V}_0 , {\bf V}_1 \rangle$ respectively. 

By modular the subgroup ${\goth G}_N$ of ${\goth G}$ in (\req(CPNg)), we obtain three normal subgroups of the CP group $G_N$ from (\req(uHrl)), with the sets of generators:
\bea(lll)
H^\circ = \{ {\sf V}_0,{\sf V}_1 \} \subset H = \langle {\sf V}_0 ,{\sf V}_1 \rangle ,& H_l^\circ = \{ {\sf V}_0,{\sf V}_2 \} \subset
H_l=\langle {\sf V}_0,{\sf V}_2 \rangle , \\
H_r^\circ = \{ {\sf V}_1,{\sf V}_2 \} \subset H_r =\langle {\sf V}_1,{\sf V}_2 \rangle.
\elea(Hrl)
Then Lemma \ref{lem:Uij} and the relation (\req(u01IJ)) yield the following results:
\begin{prop}\label{prop:Dih} (i) The CP group $G_N$ is characterized as the group with generators ${\sf V}_0, {\sf V}_1, U$ and ${\sf i}, {\sf j}$  satisfying the relations:
\bea(lllll)
{\sf V}_0 {\sf V}_1 = {\sf V}_1{\sf V}_0, & C_U ({\sf V}_k) = {\sf V}_k, & {\sf C}_{{\sf j}}({\sf V}_k) = {\sf V}_k^{-1} &C_{{\sf i}}({\sf V}_k) = {\sf V}_k^{(-1)^k}, &  (k=0, 1),  \\
 {\sf j}^2 = {\sf i}^2  = 1, &U {\sf j} = {\sf V}_0 {\sf V}_1 {\sf j}U, & U{\sf i} = {\sf V}_0 {\sf i} U^{-1}, & {\sf j} {\sf i} = {\sf V}_1 {\sf i} {\sf j}, & {\sf V}_0^N = {\sf V}_1^N = U^N = 1 ,
\elea(VUij)
where ${\sf V}_0, {\sf V}_1, U$ and ${\sf i}, {\sf j}$ are related to $U, {\sf u}_1, {\sf u}_2$ in (\req(CPNg)) by
\bea(lllll)
{\sf V}_0= U^2 {\sf u}_2^{-2}{\sf u}_1^{-2}, ~ {\sf V}_1= {\sf u}_1^2, & {\sf j} = U {\sf u}_2^{-1}, ~ ~ {\sf i} = U {\sf u}_2^{-1} {\sf u}_1, &
\Longleftrightarrow & {\sf u}_1 = {\sf j}{\sf i}, ~ {\sf u}_2 = {\sf j}U.
\elea(GN0lij)

(ii) $G_N$ is related to $\ZZ_2 \times D_N$  by
\bea(lll)
H ~ (= \langle {\sf V}_0 , {\sf V}_1 \rangle  \simeq \ZZ_N^2 )  \lhd & G_N \longrightarrow  G_N/H = \ZZ_2 \times D_N & (=  \langle \sigma \rangle \times \langle \theta, \iota \rangle) , \\
\sigma ^2= \iota^2= \theta^N =1, & \sigma \theta = \theta \sigma , ~ ~ \sigma \iota= \iota \sigma , &
\theta \iota = \iota \theta ^{-1},
\elea(GNDih)
where $\theta , \sigma, \iota$ are the classes of $U,{\sf j},{\sf i}$ in $G_N/\langle {\sf V}_0 , {\sf V}_1 \rangle$ respectively.
\end{prop}
$\Box$ \par \noindent
Since $G_N$ is a normal subgroup of $\widetilde{G}_N$, the conjugation of an element ${\sf V} \in \widetilde{G}_N$ on (\req(GN0lij)) produces a representation of (\req(VUij)). In particular, by using (\req(CPG)) and (\req(STG12)), the conjugation $C_{\sf v}$ for ${\sf v} \in G_N$ gives rise to the following representations of (\req(VUij)):
\bea(lllll)
{\sf v} \in G_N: &C_{\sf v}({\sf V}_0 , {\sf V}_1 ), & C_{\sf v}(U), & C_{\sf v}({\sf j}), &C_{\sf v}({\sf i}) ,   \\
 {\sf u}_1 :& ({\sf V}_0^{-1} , {\sf V}_1 ), & U^{-1}{\sf u}_1^2 (= {\sf V}_1 U^{-1}), & U{\sf u}_2^{-1}{\sf u}_1^{-2}  (= {\sf V}_1 ~ {\sf j}) , &U{\sf u}_2^{-1}{\sf u}_1^{-1}= ({\sf V}_1 ~ {\sf i}) ; \\
{\sf u}_1^2 :& ({\sf V}_0 , {\sf V}_1 ), & U, & U{\sf u}_2^{-1}{\sf u}_1^{-4} (= {\sf V}_1^2 ~{\sf j}), & U{\sf u}_2^{-1}{\sf u}_1^{-5}(= {\sf V}_1^2 ~ {\sf i}); \\
{\sf u}_2 :& ({\sf V}_0^{-1} , {\sf V}_1^{-1} ), & U^{-1}{\sf u}_2^2 (= {\sf V}_2 U^{-1}), &  U^{-1}{\sf u}_2
(= {\sf V}_2 U^{-2} {\sf j}),  & U{\sf u}_2 {\sf u}_1 (= {\sf V}_2 {\sf i}) ;\\
{\sf u}_2^2 :& ({\sf V}_0 , {\sf V}_1 ), & U, &  U{\sf u}_2^{-1} (= {\sf j}),  &U {\sf u}_2^3 {\sf u}_1 (= {\sf V}_2^2 {\sf i}) ;\\
U :& ({\sf V}_0 , {\sf V}_1 ), & U, &  U^3 {\sf u}_2^{-3} (= {\sf V}_0{\sf V}_1 {\sf j}), & U{\sf u}_2 {\sf u}_1^{-1}(={\sf V}_1{\sf V}_2 {\sf i}). \\
\elea(GNDi)
where ${\sf V}_0 , {\sf V}_1, {\sf j}, {\sf i}$ are in (\req(GN0lij)), and ${\sf V}_2 = {\sf u}_2^2 (= U^2{\sf V}_0^{-1} {\sf V}_1^{-1})$ in (\req(G12)). Similarly, by using (\req(ST)), (\req(tCPNrel)) and (\req(STG12)), there are two other representations of (\req(VUij)):
\bea(lllll)
{\sf V}: &C_{\sf V}({\sf V}_0 , {\sf V}_1 ), & C_{\sf V}(U), & C_{\sf V}({\sf j}), &C_{\sf V}({\sf i}) ,   \\
 S^{-1}: &({\sf V}_0 , {\sf V}_2), &U,&U{\sf u}_1^{-1},  & U{\sf u}_2^{-1} {\sf u}_1^{-1};  \\
ST^{-1}S^{-1}: &({\sf V}_2, {\sf V}_1 ), &U, &U{\sf u}_2^{-1}{\sf u}_1,  & U{\sf u}_2^{-1} {\sf u}_1^2,
\elea(GNlr)
by which one finds two $(\ZZ_2 \times D_N)$-structures related to $G_N$:
\bea(lll)
H_l  \lhd  G_N \longrightarrow  G_N/H_l = \ZZ_2 \times D_N (=  \langle \sigma_l \rangle \times \langle \theta_l, \iota_l \rangle), ~ 
(U, U{\sf u}_1^{-1}, U{\sf u}_2^{-1} {\sf u}_1^{-1}) \mapsto (\theta_l, \sigma_l, \iota_l) ; \\
H_r \lhd G_N \longrightarrow G_N/H_r=\ZZ_2\times D_N (=\langle \sigma_r\rangle \times \langle \theta_r, \iota_r \rangle),~(U, U{\sf u}_2^{-1}{\sf u}_1, U{\sf u}_2^{-1}{\sf u}_1^2) \mapsto (\theta_r, \sigma_r, \iota_r). 
\elea(GNlrDi)

We now study the change of structures in (\req(GNDih)) under the conjugation action of $\langle S, T \rangle (\subseteq \widetilde{G}_N)$ in (\req(conjtGN)). 
\begin{lem}\label{lem:HsST} (i) For $N=2$, the subgroups in (\req(Hrl)) are equal: $H= H_l = H_r$, which is equal to the normal subgroup $G'_{2, 1}$ of $\widetilde{G}_2$ in (\req(GN12)).

(ii) For $N \geq 3$, the subgroups in (\req(Hrl)) are all distinct. Furthermore, $C_{\sf V}(H)= H_l$ (or $H_r$) if and only if $C_{\sf V}(H^\circ)= H_l^\circ$ (or $H_r^\circ$ respectively), where ${\sf V}$ is an element in $\widetilde{G}_N$.
\end{lem}
{\it Proof.} Note that the equality of a pair subgroups in (\req(Hrl)) is equivalent to $H= H_l = H_r= G'_{N, 1}$. 
When $N=2$, $U^2=1$, hence follows $(i)$. When $N \geq 3$, the order of $G'_{N, 1} > N^2$, by which the three subgroups in (\req(Hrl)) are distinct. If $C_M (H)= H_l$ and $C_M (H^\circ)\neq H_l^\circ$ for some $M \in \langle S, T \rangle$, using (\req(STG12)),  one finds $C_M (H)=  G'_{N, 1}$, which contradicts $N \geq 3$. Hence follows $(ii)$.
$\Box$ \par \vspace{.2in} \noindent
Using Lemma \ref{lem:HsST}, we can determine the normalizer of $H, H_l, H_r$ in $\widetilde{G}_N$ for $N \geq 3$.
\begin{prop}\label{prop:MHlr} For ${\sf V} \in \widetilde{G}_N$, let $\overline{\sf V}$ be the element in $SL_2(\ZZ_2)$  corresponding to the class of  ${\sf V}$ in (\req(GNSL22)). Then 
\bea(lll)
C_{\sf V}(H^\circ) = H^\circ & \leftrightarrow &\overline{\sf V} = 1,  \overline{T} , \\
C_{\sf V}(H^\circ) = H_l^\circ & \leftrightarrow &\overline{\sf V} = \overline{S}, \overline{ST}, \\
C_{\sf V}(H^\circ) = H_r^\circ & \leftrightarrow &\overline{\sf V} = \overline{STS^{-1}}, \overline{STS^{-1}T},  \\
\elea(HHrlr)
where $\overline{S}, \overline{T}$ are the generators of $SL_2(\ZZ_2)$ in (\req(ST-)) with $\overline{S}^2=\overline{T}^2=1$.
As a consequence, for $N \geq 3$, the equivalent relations (\req(HHrlr)) are still valid when replacing $H^\circ, H_l^\circ, H_r^\circ$ to 
$H, H_l, H_r$.
\end{prop}
{\it Proof.} By (\req(STG12)), $H^\circ, H_l^\circ, H_r^\circ$ are invariant under the conjugation of an element in $\langle G_N, T^2, ST^{-2}S^{-1} \rangle$. Hence the transformation relations on the left hand side of (\req(HHrlr)) depend only on the quotient group $\widetilde{G}_N/\langle G_N, T^2, ST^{-2}S^{-1} \rangle$, which is isomorphic to $SL_2(\ZZ_2)$ by (\req(GNSL22)). The relation (\req(STG12)) yields $C_T (H^\circ) = H^\circ$, $C_S (H^\circ) = H_l^\circ$, $C_{STS^{-1}}(H^\circ)$, hence follows the conclusion.   
$\Box$ \par \vspace{.1in} \noindent
\begin{cor} \label{cor:nHH} For $N \geq 3$, the normalizer of $H$ in $\widetilde{G}_N$, and $N(H; H_l), N(H; H_l)$ in (\req(BNH)), are  given by
\bea(lll)
N(H) = \langle G_N, T^{-1}, ST^2 S^{-1} \rangle, &
N(H; H_l) = S^{-1} N(H) , &N(H; H_r) = ST^{-1}S^{-1} N(H). 
\elea(NH)
By conjugation of $S^{-1}, ST^{-1}S^{-1}$ on $N(H)$, one obtains the normalizers of $H_l, H_r$:
$$
N(H_l) = \langle G_N, ST^{-1}S^{-1}, T^2 \rangle, ~ N(H) = \langle G_N, T^2S, ST^2 S^{-1} \rangle.
$$
\end{cor}
$\Box$ \par \vspace{.1in} \noindent
We now examine the relationship of $(\ZZ_2 \times D_N)$-structure for $G_N/H$ in (\req(GNDih)) under the conjugation of $N(H)$ in (\req(NH)). The change of $\theta, \sigma, \iota$ in (\req(GNDih)) depends on the representation of (\req(VUij)) and (\req(GN0lij)) under $N(H)$-conjugation, in which the $G_N$-conjugation are determined by the representations in (\req(GNDi)). By (\req(NH)), it remains to examine the change of (\req(GN0lij)) under $\langle T^{-1}, ST^2 S^{-1} \rangle$-conjugation. 
\begin{lem}\label{lem:NHDih} (i) Let $M$ be an element in $\langle T^{-1}, ST^2 S^{-1} \rangle$. The change of (\req(GN0lij)) under the conjugation $C_M$ of $M$ is given by
$$
C_M({\sf V}_0, {\sf V}_1) = \left\{\begin{array}{ll}({\sf V}_0, {\sf V}_1) & {\rm if} ~ M \in \langle T^{-2}, ST^2 S^{-1} \rangle, \\
({\sf V}_1, {\sf V}_0) & {\rm if} ~ M = T^{2n+1} ,
 \end{array}\right.
$$
and 
\bea(lllll)
M: & C_M(U), & C_M({\sf j}), &C_M({\sf i}) ,   \\
T^{-k}: &U,& {\sf j},  & {\sf u}_2^{-k} {\sf i} ~ ;  \\
ST^{2k}S^{-1}:&U, &{\sf u}_1^{2k} {\sf j} ,  & {\sf u}_1^{2k} {\sf i}  .
\elea(NHUij)
where ${\sf V}_0 , {\sf V}_1, {\sf j}, {\sf i}$ are defined in (\req(GN0lij)). 

(ii) The following equalities hold for elements in $\langle T^{-1}, ST^2 S^{-1} \rangle$ and $G_N$: 
\bea(lll)
T^{-4} = {\sf u}_2^{-2}, ~  ST^4 S^{-1}= {\sf u}_1^2, & (ST^2S^{-1})T^{-1}= T^{-3}(ST^2S^{-1}) U{\sf u}_1^{-2}, \\
(ST^2 S^{-1}) T^{-2} = T^{-2} (ST^2 S^{-1}) , &
(ST^2S^{-1})T^{-3}= T^{-1}(ST^2S^{-1}) U{\sf u}_1^{-2}{\sf u}_2^{-2} 
\elea(NHGN)
by which $G_N \cap \langle T^{-1}, ST^2 S^{-1} \rangle = G_{N, 1}$, and the quotient group $G_N \langle T^{-1}, ST^2 S^{-1} \rangle/ G_N$ consists of 8 elements represented by $T^{-k}, ST^2 S^{-1}T^{-k}$. 
\end{lem}
{\it Proof.} It is easy to see that $(i)$ follows from (\req(ST)), (\req(tCPNrel)) and (\req(STG12)). By (\req(Gt/GN)) and the first relation in (\req(HHrlr)), $|G_N \langle T^{-1}, ST^2 S^{-1} \rangle/ G_N| = 8$. On the other hand, 
by (\req(UGtAut)), (\req(uCPrel)), the computation via formulas in (\req(ST)) implies  equalities in (\req(NHGN)) true for the universal CP group $\widetilde{\goth G}$, hence also valid for $\widetilde{G}_N$. Then follows $(ii)$.
$\Box$ \par \vspace{.1in} \noindent
{\bf Remark.} Note that representations of $G_N$ in (\req(GNlr)) yield the $(\ZZ_2 \times D_N)$-structures of $G_N/H_l, G_N/H_r$, hence a similar relation of $(\ZZ_2 \times D_N)$-structure change under $N(H_l), N(H_r)$-conjugation can be obtained from Lemma \ref{lem:NHDih} by applying the conjugation of $S^{-1}, ST^{-1}S^{-1}$.
$\Box$ \par \vspace{.2in} \noindent

\section{Fermat Surface and the Family of Rapidity Curves in Chiral Potts Model \label{sec.Fermat}}
In this section, we make a thorough investigation about the algebraic geometry structure of the rapidity family (\req(WFaml)) and its related family of hyperelliptic curves in $\tau^{(2)}$-model. The symmetry of these fibrations is studied through the geometrical representation of the structure group introduced in Section \ref{sec:Alg}.

\subsection{Fermat surface and symmetry group of the rapidity family  in Chiral Potts Model \label{subs.FSym}}
\setcounter{equation}{0}
The family in ${\goth W}$ in (\req(WFaml)) is a subvariety in $\Lambda \times \PZ^3$ with a projection to $\PZ^3$,
$$
{\goth W} \subset \Lambda \times \PZ^3  \longrightarrow \PZ^3 .
$$
Indeed, the above projection induces an isomorphism between ${\goth W}$ and a Fermat hypersurface of $\PZ^3$: 
\begin{prop}\label{prop:Ferm} The complete family (\req(WFaml)) of CP rapidity curves is isomorphic to the following Fermat surface  in $\PZ^3$:
\bea(lll)
{\goth W} ~ (= {\goth W}^{(N)}) &\simeq &\{ [a, b, c, d] \in \PZ^3 | a^{2N} + c^{2N} = b^{2N} + d^{2N} \} \\
\downarrow & & \pi \downarrow \\
\Lambda & =  & ~ ~ \Lambda .
\elea(WFerm)
\end{prop}
{\it Proof.} For convenience, in this proof, we denote he Fermat surface on the right of (\req(WFerm)) by ${\goth F}$. It is obvious that ${\goth W}_{k', k} \subset {\goth F}$ for all $(k', k) \in \Lambda$. Indeed,  when $k \neq 0$, ${\goth W}_{k', k}$ in (\req(Wk'k)) or (\req(Wdeg)) is defined by the right two equations of (\req(Wk'k)) in the form:  
\bea(llll)
a^N = - \frac{k'}{k}c^N + \frac{1}{k} d^N , &
b^N =  \frac{1}{k} c^N -  \frac{k'}{k} d^N ,
\elea(WinF) 
by which $a^{2N}- b^{2N} = \frac{1- k'^2}{k^2} (-c^{2N}+d^{2N})$. The constraint of $(k', k)$ for $\Lambda$ in (\req(kk'kap)) is equivalent to the Fermat relation of ${\goth F}$. When $k = 0$,  ${\goth W}_{\pm 1, 0}$ in (\req(Wdeg)) consists of rational curves in ${\goth F}$ defined by $(a^N, c^N) = \pm (-b^N, d^N)$ respectively. In order to show the one-to-one correspondence in (\req(WFerm)), it suffices to construct the projection $\pi$ from ${\goth F}$ to $\Lambda$ in (\req(kk'kap)). For an element $p= [a, b, c, d] \in {\goth F}$, one can find $k', k \in \PZ^1$ satisfying the relation (\req(WinF)) when $c^{2N} - d^{2N} \neq 0$, equivalently,  $(\frac{-k'}{k}, \frac{1}{k}) = \frac{1}{c^{2N} - d^{2N}}(c^N a^N- d^N b^N, -d^Na^N, c^N b^N)$, hence $\pi (p):= (k',  k) \in \Lambda$. When $c^{2N} - d^{2N} = 0$, then $a^{2N} - b^{2N} = 0$, equivalent to $(a^N, c^N) = \pm (-b^N, d^N)$ or $\pm (b^N, d^N)$, where in the formal case, we define  $\pi (p) := (\pm 1, 0)$. When $(a^N, c^N) = \pm (b^N, d^N) $ but $\neq \pm (-b^N, d^N)$, then $a^N \neq 0, d^N \neq 0$, where $\pi (p)= (k', k) \in \Lambda$ is defined by the relation $(\frac{1\mp k'}{1 \pm k'}, \frac{1 \mp k'}{k})= (\frac{a^{2N}}{d^{2N}}, \frac{a^N}{d^N})$. Furthermore, by the construction of $\pi$, we find the fiber $\pi^{-1}(k',  k)$ in ${\goth F}$ is equal to ${\goth W}_{k', k}$. 
$\Box$ \par \noindent  
{\bf Remark.} In the proof of the above Proposition \ref{prop:Ferm}, we find that $(a^N, c^N) = \pm (b^N, d^N)$ define the $2N^2$ "horizontal" lines of the fibration (\req(WFerm)), whose intersection with ${\goth W}_{k, k}$ in (\req(Wk'k)) are the vertical rapidity in superintegrable CPM \cite{AMP, B88, GR, R1003}. 
$\Box$ \par \noindent \vspace{.1in}
 From now on, we shall identify the family of CP rapidity curves, ${\goth W}$ in (\req(WFaml)), with the Fermat hypersurface in (\req(WFerm)). We now identify the modular ($N$-state) CP group ${\goth G}_N $ in (\req(CPNmg)) with the automorphism group of the fibration ${\goth W}$ over $\Lambda$ in (\req(WFerm)).  
Represent the generators of $\widetilde{G}_N$ in (\req(CPNmg)) by the following automorphisms of ${\goth W}$: 
\bea(ll)
{\sf u}_1: [a , b , c, d] \mapsto [\omega^\frac{1}{2} d, c , b, \omega^\frac{1}{2}  a] , &
{\sf u}_2: [a , b , c, d]  \mapsto [b , \omega a , d, c] , \\
U : [a , b , c, d]  \mapsto [\omega a , b , c, d] ,  \\
S: [a, b, c, d]  \mapsto [\omega^\frac{-1}{2} a, d, c, b ], &
T: [a, b, c, d]  \mapsto [\omega^\frac{1}{4} a, \omega^\frac{1}{4} b, d, c], 
\elea(GtWf)
where $\omega^\frac{1}{2} = e^\frac{\pi {\rm i}}{N}$, $\omega^\frac{1}{4}= e^\frac{\pi {\rm i}}{2N}$. One finds that
the relations (\req(CPG)) (\req(CPN)) and (\req(tCPNrel)) hold. Hence we obtain a representation of the modular ($N$-state) CP group $\widetilde{G}_N$ as an automorphism group of ${\goth W}$ . The generators of the subgroups $G_N$ in (\req(CPNg)) or $G'_{N, 1}$ in (\req(GN12)) of ${\goth G}_N$ are related to automorphisms in (\req(Sym)) of a CP rapidity curve ${\goth W}_{k', k}$ in (\req(Wk'k)) or (\req(Wdeg)): 
\bea(lll)
M^{(1)}={\sf u}_1^2 , ~ M^{(2)}= {\sf u}_2^2, & M^{(3)}= U {\sf u}_2 {\sf u}_1, &  M^{(4)}= U^{-1}{\sf u}_2^2 {\sf u}_1^2 , ~ ~ M^{(5)}= U^{-1} {\sf u}_2^2 {\sf u}_1 ; \\
{\sf u}_1 = {M^{(5)}}^{-1} M^{(4)}, & {\sf u}_2 =   M^{(3)}M^{(5)}, & U =M^{(1)}M^{(2)}{M^{(4)}}^{-1} ; \\
{\sf V}_0= {M^{(0)}}, & {\sf V}_1 = M^{(1)}, & {\sf V}_2 = M^{(2)}.
\elea(uVMaut) 
Hence by (\req(AutWk)), $G_N$ can be identified with the automorphism group of ${\goth W}_{k', k}$ in (\req(Wk'k)) when $N \geq 3$:
\bea(lll)
G_N = \langle {\sf u}_1, {\sf u}_2, U \rangle = {\rm Aut}({\goth W}_{k', k}), & k' \neq 0, \pm 1, \pm \infty, & (N \geq 3).
\elea(GtNWk)
By (\req(slovN)), ${\rm Aut}({\goth W}_{k', k})$ is a solvable group\footnote{The solvable group structure of ${\rm Aut}({\goth W}_{k', k})$ in (\req(slovN)) here is different from that in \cite{R04} Proposition 1.} of order $4N^3$. For $N \geq 2$, the transformations $S, T$ in (\req(GtWf)) are  automorphisms  of the fibration ${\goth W}$ over $\Lambda$  in (\req(WFerm)) with the properties
\bea(ll)
S: {\goth W}_{k', k} \simeq {\goth W}_{k, k'},  & T: {\goth W}_{k', k} \simeq {\goth W}_{\frac{1}{k'}, \frac{{\rm i} k}{k'}}.
\elea(STWk')
Indeed, the quotient group in (\req(Gt/GN)) 
$$
\widetilde{G}_N/G_N = \langle \overline{S}, \overline{T} \rangle = PSL_2(\ZZ_4)
$$ 
acts on the parameter space $\Lambda$, where $\overline{S}, \overline{T}$ are the classes of $S, T$ in $\widetilde{G}_N/G_N$, identified with the generators of $PSL_2(\ZZ_4)$ in (\req(STSL4)). With the identification of $(k', k) \in \Lambda$ and $\kappa' \in \PZ^1$ in (\req(kk'kap)), the action of $\widetilde{G}_N/G_N$ on $\Lambda$ are generated by
\bea(lll)
\overline{S} : (k', k) \mapsto (k, k') ~ ~ ( \Leftrightarrow  \kappa' \mapsto {\rm i} \kappa'^{-1} ), & 
\overline{T} : (k', k) \mapsto (\frac{1}{k'}, \frac{{\rm i} k}{k'})~ ~ (\Leftrightarrow  \kappa' \mapsto \frac{1+ {\rm i} \kappa'}{\kappa' + {\rm i}}).
\elea(STlambda) 
\begin{lem}\label{lem:PSL2} The transformations in (\req(STlambda)) give rise to the identification of $\widetilde{G}_N/G_N$ with the automorphism group of $\Lambda ~ (\simeq \PZ^1) $ in (\req(kk'kap)) preserving the degenerated-parameter set  $\Lambda_{\rm deg}$ in (\req(Wdeg)):
\bea(lll)
\widetilde{G}_N/G_N = PSL_2(\ZZ_4) = {\rm Aut}(\Lambda,  \Lambda_{\rm deg}) (:= \{ \phi \in {\rm Aut}(\Lambda ) | \phi (\Lambda_{\rm deg}) = \Lambda_{\rm deg}  \} ),
\elea(PSL2)
where $\Lambda_{\rm deg} := \{(k', k) \in \Lambda | k'= \infty, 0, \pm 1 \} ~ (\simeq  \{ \kappa' =  0, \infty, \pm 1, \pm {\rm i} \in \PZ^1  \})$.
\end{lem}
{\it Proof.} With the identification in (\req(kk'kap)): $\Lambda = \PZ^1 $,  ${\rm Aut}(\Lambda,  \Lambda_{\rm deg}) \subset {\rm Aut}(\PZ^1)$. Note that there are three types of degenerated parameters: 
$$
\begin{array}{lll}
\Lambda_{\rm deg} =  \sqcup_{k=\infty, \pm 1} \Lambda_{{\rm deg}, k} & = \Lambda_{{\rm deg} , \infty} \sqcup \Lambda_{{\rm deg} , 1} \sqcup \Lambda_{{\rm deg} , -1} , \\
\Lambda_{{\rm deg} , \infty} :=\{ \kappa' = 0^{\pm 1} \}, & \Lambda_{{\rm deg} , 1} := \{ \kappa' = \pm 1 \}, &\Lambda_{{\rm deg} , -1} := \{ \kappa' = \pm {\rm i} \}.
\end{array}
$$
Consider the following subgroup of ${\rm Aut}(\Lambda,  \Lambda_{\rm deg})$:
$$
\begin{array}{lll}
{\rm Aut}(\Lambda,  \sqcup_k \Lambda_{{\rm deg}, k}) := \{ \phi \in {\rm Aut}(\Lambda ) | \phi (\Lambda_{{\rm deg}, k}) = \Lambda_{{\rm deg}, \sigma(k)}, ~ \sigma: {\rm a ~ permutation ~ of ~} k' (= \infty, \pm 1) \}. 
\end{array}
$$
By (\req(PSL2)), one finds $\overline{S}, \overline{T} \in {\rm Aut}(\Lambda,  \sqcup_k \Lambda_{{\rm deg}, k})$: 
\bea(lll)
\overline{S} (0^{\pm 1} , \pm 1, \pm {\rm i}) = (0^{\mp 1} , \pm {\rm i}, \pm 1 ) , & 
\overline{T} (0^{\pm 1} , \pm 1, \pm {\rm i}) = ( \mp {\rm i}, \pm 1,  0^{\pm 1}) ; & \overline{S}^2= \overline{T}^4 = 1; 
\elea(Psl24)
and the normal subgroup $\langle \overline{T}^2,  \overline{S} \overline{T}^2\overline{S} \rangle (\simeq \ZZ_2^2 ) $ of $PSL_2(\ZZ_4)$ is expressed by
$$
\overline{T}^2 : \kappa' \mapsto \kappa'^{-1} ; ~  \overline{S} \overline{T}^2\overline{S}: \kappa' \mapsto -\kappa'^{-1} ; ~ (\overline{S} \overline{T}^2)^2: \kappa' \mapsto -\kappa'
$$
in ${\rm Aut}(\Lambda,  \sqcup_k \Lambda_{{\rm deg}, k})$. Indeed, one finds
$$
\begin{array}{lll}
\langle \overline{T}^2,  \overline{S} \overline{T}^2\overline{S} \rangle = \{ \phi \in {\rm Aut}(\Lambda,  \sqcup_k \Lambda_{{\rm deg}, k})| \phi (\Lambda_{{\rm deg}, k}) = \Lambda_{{\rm deg}, k}  \} \lhd {\rm Aut}(\Lambda,  \sqcup_k \Lambda_{{\rm deg}, k}) 
\end{array}
$$
by which through the morphism $\PZ^1 \longrightarrow \PZ^1 , \kappa' \mapsto \frac{\kappa'^2+ \kappa'^{-2}}{2}$, the quotient group can be embedded into the automorphism group of $\PZ^1$ permuting three elements $ \infty, \pm 1$:
$$
\begin{array}{lll}
{\rm Aut}(\Lambda,  \sqcup_k \Lambda_{{\rm deg}, k})/\langle \overline{T}^2,  \overline{S} \overline{T}^2\overline{S} \rangle \hookrightarrow {\rm Aut}(\PZ^1, \{ \infty, \pm 1 \}).  
\end{array}
$$
The relation (\req(Psl24)) implies that $\overline{S}, \overline{T}$ induce the automorphism $[\overline{S}], [\overline{T}]$ in ${\rm Aut}(\PZ^1, \{ \infty, \pm 1 \})$ with $[\overline{S}](\infty, \pm 1 )= (\infty, \mp 1)$ and 
$[\overline{T}](\infty, 1, -1 )= (-1, 1, \infty)$. Since ${\rm Aut}(\PZ^1, \{ \infty, \pm 1 \})$ is isomorphic to the permutation group of $\{ \infty, \pm 1 \}$, we obtain
$$
SL_2(\ZZ_2) \simeq PSL_2(\ZZ_4)/\langle \overline{T}^2,  \overline{S} \overline{T}^2\overline{S} \rangle \simeq {\rm Aut}(\Lambda,  \sqcup_k \Lambda_{{\rm deg}, k})/\langle \overline{T}^2,  \overline{S} \overline{T}^2\overline{S} \rangle \simeq  {\rm Aut}(\PZ^1, \{ \infty, \pm 1 \}),
$$
hence $PSL_2(\ZZ_4) = {\rm Aut}(\Lambda,  \sqcup_k \Lambda_{{\rm deg}, k})$. It remains to show 
$$
{\rm Aut}(\Lambda,  \sqcup_k \Lambda_{{\rm deg}, k})= {\rm Aut}(\Lambda,  \Lambda_{\rm deg}).
$$ 
Otherwise, there is an element $\phi \in {\rm Aut}(\Lambda,  \Lambda_{\rm deg}) \setminus {\rm Aut}(\Lambda,  \sqcup_k \Lambda_{{\rm deg}, k})$. By composing with some elements in ${\rm Aut}(\Lambda,  \sqcup_k \Lambda_{{\rm deg}, k})$, we may assume $\phi (\infty ) = \infty, \phi (0) = 1$. Then $\phi (\kappa') = \alpha \kappa' + 1$ for some $\alpha \in \CZ$ with $\{ \pm \alpha + 1, \pm{\rm i} \alpha + 1 \} =\{ 0, -1, \pm {\rm i} \}$, which leads to a contradiction.
$\Box$ \par \noindent \vspace{.1in} 
{\bf Remark.} There are 24 elements $\phi \in {\rm Aut}(\Lambda,  \Lambda_{\rm deg})$ with $\phi (k', k) $ of $(k', k) \in \Lambda$  given by
\bea(lll)
(k', \pm k), (\frac{1}{k'}, \pm \frac{{\rm i} k}{k'}), (-k', \pm k), (\frac{-1}{k'}, \pm \frac{{\rm i} k}{k'}),  (\Longleftrightarrow ~ \kappa'^{\pm 1}, (\frac{1 + {\rm i} \kappa'}{\kappa' + {\rm i}})^{\pm 1}, - \kappa'^{\mp 1}, (\frac{\kappa'-{\rm i}}{{\rm i}\kappa' -1 })^{\mp 1})  \\
(k, \pm k'), (\frac{1}{k}, \pm \frac{{\rm i} k'}{k}), (-k, \pm k'), (\frac{-1}{k}, \pm \frac{{\rm i} k'}{k}), ( \Longleftrightarrow ~ \pm {\rm i} \kappa'^{\mp 1 },  \pm {\rm i}(\frac{1-\kappa'}{1+\kappa'})^{\pm 1}, \pm {\rm i}\kappa'^{\pm 1},  \pm {\rm i} (\frac{1+\kappa'}{1-\kappa'})^{\pm 1})  \\
(\frac{{\rm i} k}{k'}, \pm \frac{1}{k'}), (\frac{-{\rm i} k}{k'}, \pm \frac{1}{k'}), 
(\frac{{\rm i} k'}{k}, \pm \frac{1}{k}), (\frac{-{\rm i} k'}{k}, \pm \frac{1}{k}),  ( \Longleftrightarrow ~ 
(\frac{{\rm i} \kappa'-1}{{\rm i} \kappa' + 1})^{\pm 1}, (\frac{1-{\rm i} \kappa'}{1+{\rm i} \kappa'})^{\mp 1}, 
(\frac{1+  \kappa'}{1- \kappa'})^{\pm 1}, (\frac{\kappa'+1}{\kappa'-1})^{\mp 1}). 
\elea(24k')
Write $\phi = \overline{M}$ for some $M \in \langle S, T \rangle$ in (\req(PSL2)), then $M$ gives rise to an isomorphism of the fibration ${\goth W}$ over $\Lambda$ with  $M : {\goth W}_{k', k} \simeq {\goth W}_{\overline{M} (k, k')}$ as those in (\req(STWk')).
$\Box$ \par \noindent \vspace{.1in} 
Using Lemma \ref{lem:PSL2}, we now show
\begin{thm}\label{thm:WGtN} For $N \geq 2$, the modular CP group $\widetilde{G}_N$ in the representation (\req(GtWf)) is identified with the automorphism group of the fibration ${\goth W}^{(N)}$ over $\Lambda$ in (\req(WFerm)) with the order $96 N^3$.
\end{thm}
{\it Proof.} It is obvious that the representation (\req(GtWf)) embeds $\widetilde{G}_N$ into the automorphism subgroup of the fibration ${\goth W}^{(N)}$ over $\Lambda$. Suppose $\Phi$ is an automorphism of ${\goth W}^{(N)}$ over $\Lambda$. We are going to show $\Phi \in \widetilde{G}_N$. Since $\Phi$ induces an automorphism of $\Lambda$ preserving the degenerated-parameter set $\Lambda_{\rm deg}$, by (\req(PSL2)) and composing with some automorphism in $\widetilde{G}_N$, we may assume $\Phi$ induces the identity on the base $\Lambda$. Hence $\Phi$ preserves each fiber ${\goth W}_{k', k}$ of ${\goth W}$. When $N \geq 3$, the relation (\req(GtNWk)) yields $\Phi \in G_N$. It remains the case $N=2$, where ${\goth W}^{(2)}$ is a K3 surface with an elliptic fibration over $\Lambda$. It is known that line bundles of ${\goth W}^{(2)}$ are described by cohomology elements in
$$
H^1({\goth W}^{(2)}, {\cal O}^*) \simeq H^2({\goth W}^{(2)}, \ZZ ) = \ZZ^{20}.
$$
The hyperplane section of $\PZ^3$ and a general fiber of the fibration (\req(WFerm)) give rise to two cohomology elements of ${\goth W}^{(2)}$, denoted by $[h], [f]$ respectively. The rest cohomology ($\QZ$-)basis elements are contributed from the 
six degenerated fibers, ${\goth W}^{(2)}_{\kappa'} ~ (\kappa' \in \Lambda_{\rm deg})$ in (\req(Wdeg)). Each consists of four lines, $\ell^{(\pm , \pm )}$,  intersecting normally only at $\ell^{(+, +)} \cdot \ell^{(+, -)}= \ell^{(+, +)} \cdot \ell^{(-, +)} = \ell^{(-, +)} \cdot \ell^{(-, -)} = \ell^{(+, -)}\cdot \ell^{(-, -)} =1$, and the four lines $\ell^{(\pm , \pm )}$'s give rise to four basis elements subject to the cohomologous equivalent relation: $\sum \ell^{(\pm , \pm )} \sim [f]$. The morphism  of cohomologous group induced by $\Phi$ leaves $[f]$ invariant, and permutes the four lines in each degenerated fiber, hence leaves $[h]$ invariant. As a consequence, the automorphism $\Phi$ of ${\goth W}^{(2)}$ is induced from a projective linear transformation of $\PZ^3$, which preserves relations in (\req(Wk'k)) and (\req(Wdeg)) for all $(k', k) \in \Lambda$. Hence $\Phi$ induces a permutation of $(a^2, b^2, c^2, d^2)$ up to some scalars. By composing with an automorphism ${\sf u}_2, U {\sf u}_2 {\sf u}_1 (=M^{(3))}$ or $U^{-1} {\sf u}_2^2 {\sf u}_1 (=M^{(5)})$ in $G_2$, we may assume $\Phi $ leaves $a^2$ invariant, by which one finds all $a^2, b^2, c^2, d^2$ invariant under $\Phi$. Hence $\Phi \in \langle {\sf u}_1^2, {\sf u}_2^2, U \rangle$.  This shows $\widetilde{G}_N$ is the automorphism group of  ${\goth W}^{(N)}$ over $\Lambda$ for all $N$, whose order is given by (\req(tGNrel)). 
$\Box$ \par \noindent \vspace{.1in} 
{\bf Remark.} In the case $N=2$, every ${\goth W}^{(2)}_{k' k}$ in (\req(Wk'k)) is an elliptic curve with infinity many symmetries. However only $32$ symmetries (in $G_2$) are preserved in the elliptic family ${\goth W}^{(2)}$, which are  contained in the symmetry group $\widetilde{G}_2$. 
In section \ref{subs.K3}, we shall discuss the elliptic K3 surface ${\goth W}^{(2)}_{k' k}$ via uniformization of elliptic curves.
$\Box$ \par \vspace{.1in} 

We now describe all the  projective lines of $\PZ^3$ in (\req(WFerm)).  
\begin{prop}\label{prop:lineF} The projective lines in the Fermat hypersurface (\req(WFerm)) and  degenerated fibers in (\req(Wdeg)) are related by 
\bea(lll)
a^{2N} - b^{2N} =  c^{2N}-d^{2N}= 0,  & ( \supset {\goth W}_{\pm 1, 0} ), \\
a^{2N} - d^{2N}  = b^{2N} - c^{2N} = 0 & ( \supset {\goth W}_{0, \pm 1} ) , \\
a^{2N} + c^{2N} = b^{2N} + d^{2N} = 0 & ( \supset {\goth W}_{\infty, \infty_{\pm} }),  \\
\elea(lines)
where the zero-locus of each set of equations consists of $4N^2$ lines permuted by $G_N$ with three $G_N$-orbits. 
One-half lines form two $G_N$-orbits as the degenerated fibers in (\req(lines)); the other $G_N$-orbit consists of the rest $2N^2$ lines as the horizontal lines of the fibration (\req(Wdeg)), each of which is a $N$-fold cover under $\pi$  over $\Lambda$ branched (only) at the two degenerated parameters related to the set. Furthermore, the singular set of the degenerated fibers are given by
\bea(lll)
{\rm Sing}( {\goth W}_{\pm 1, 0}) = & {\goth W}_{\pm 1, 0} \cap \{ a=b=0 ~ {\rm or} ~ c=d=0 \} ; \\
{\rm Sing}( {\goth W}_{0, \pm 1}) = & {\goth W}_{0, \pm 1} \cap \{ a=d=0 ~ {\rm or} ~ b= c =0 \} ; \\
{\rm Sing}( {\goth W}_{\infty, \infty_{\pm}}) = & {\goth W}_{\infty, \infty_{\pm}} \cap \{ a=c=0 ~ {\rm or} ~ b= d =0 \}.
\elea(SWdeg)
\end{prop}
{\it Proof.} Note that $G_N$ leaves each set of equations in (\req(lines)) invariant, hence permutes the $4N^2$ lines of  its zero-locus, which contains two degenerated fibers as two $G_N$-orbits. First, we consider the first set of equations in (\req(lines)) with the zero-locus  
\bea(lll)
{\goth L} := \bigcup_{i, j \in \ZZ_2} {\goth L}_{i, j} , & {\rm where} ~ {\goth L}_{i, j} :& a^{N} - (-1)^i b^{N} =  c^{N}-(-1)^j d^N = 0, 
\elea(lineij)
on which $G_N$ in (\req(GtNWk)) acts by 
$$
\begin{array}{ll}
{\sf u}_1^{-1}: {\goth L}_{i, j} \longrightarrow {\goth L}_{j-1, i-1} , & {\sf u}_2^{-1}, U: {\goth L}_{i, j} \longrightarrow {\goth L}_{i, j}. 
\end{array}
$$
Hence there are three $G_N$-orbits in ${\goth L}$: ${\goth L}_{0, 1} (={\goth W}_{1, 0}), {\goth L}_{1, 0} (={\goth W}_{-1, 0})$ and ${\goth L}_{0, 0}\cup {\goth L}_{1, 1}$. One finds ${\goth L}_{i, j} \cap {\goth L}_{i', j'} \neq \emptyset$ if and only if $i \equiv i'$ or $j \equiv j' $ in $\ZZ_2$, with non-empty intersections at 
$$
\begin{array} {ll}
{\goth L}_{i, j} \cap {\goth L}_{i, j'} = \{ a^{N} - (-1)^i b^{N}= c = d = 0 \} , &{\rm when} ~  j \not\equiv j' \pmod{2}, \\
{\goth L}_{i, j} \cap {\goth L}_{i', j} = \{ a =  b = c^{N}-(-1)^j d^N  = 0 \} , & {\rm when} ~  i \not\equiv i' \pmod{2},
\end{array}
$$
the union of which is the singular set of ${\goth L}_{i, j}$. Hence follows the singular set of ${\goth W}_{\pm 1, 0}$ in (\req(SWdeg)). Indeed, each ${\goth L}_{i, j}$ consists of $N^2$ lines, 
\bea(lll)
{\goth L}_{i, j} = \bigcup_{m, n \in \ZZ_N}{\goth L}_{i, j}^{m, n}, & {\goth L}_{i, j}^{m, n}: a - \omega^{\frac{i}{2}+m} b =  c -\omega^{\frac{j}{2}+n} d = 0, ~ ~ (i, j = 0, 1) , 
\elea(Lijmn)
and ${\goth L}_{i, j}^{m, n} \cap {\goth L}_{i, j}^{m', n'} \neq \emptyset$ if and only if $m \equiv m'$ or $n \equiv n' $ with  
$$
\begin{array} {ll}
{\goth L}_{i, j}^{m, n} \cap {\goth L}_{i, j}^{m, n'} = \{ a - \omega^{\frac{i}{2}+m} b =  c = d = 0 \} , &{\rm when} ~  n \not\equiv n' \pmod{N}, \\
{\goth L}_{i, j}^{m, n} \cap {\goth L}_{i, j}^{m', n} = \{ a =  b= c -\omega^{\frac{j}{2}+n} d = 0 \} , & {\rm when} ~  m \not\equiv m' \pmod{N}.
\end{array}
$$
The $G_N$-orbit ${\goth L}_{0, 0}\cup {\goth L}_{1, 1}$ consists of $2N^2$ lines with ${\goth L}_{0, 0} \cap {\goth L}_{1, 1} = \emptyset$. For $(i, j)=(0, 0), (1, 1)$, each ${\goth L}^{m, n}_{i, j}, (m, n \in \ZZ_N)$  is a horizontal lines of (\req(Wdeg))  since the restriction of  $\pi$ defines a $N$-fold cover of ${\goth L}^{m, n}_{i, j}$ over $ \Lambda $ branched only at $(k', k)= (\pm 1, 0)$ with ${\goth L}^{m, n}_{i, j}/\langle {\sf u}_2^2 \rangle = \Lambda$, where $ {\sf u}_2$ is given in (\req(GtWf)). Therefore we obtain the conclusion for the lines in the first set of (\req(lines)). Then follow the results about lines in second and third set in (\req(lines)), which are isomorphic to those in the first set via the automorphisms $S, STS^{-1}$ in (\req(GtWf)).   
$\Box$ \par \noindent \vspace{.1in} 
{\bf Remark.} (I). By (\req(SWdeg)), there are two sets of singularities for each degenerated fiber in (\req(Wdeg)),  represented by the $a$-value $0$ or 1. For convenience, for $(k', k)= (\pm 1, 0), (0, \pm 1), (\infty, \infty_{\pm})$, we denote 
\bea(lll)
{\rm Sing}^0( {\goth W}_{k', k}) =  {\rm Sing}( {\goth W}_{k', k}) \cap \{ a =0 \} , &
{\rm Sing}^1( {\goth W}_{k', k}) = {\rm Sing}( {\goth W}_{k', k})  \cap \{ a \neq 0 \} .
\elea(SW01)
By using (\req(Lijmn)) for $(i, j)= (0, 1),(1, 0)$, and its composition with automorphisms $S, STS^{-1}$, one finds each degenerated fiber in (\req(Wdeg)) consists of $N^2$ lines:
\bea(lll)
{\goth W}_{k', k} = \bigcup_{m, n \in \ZZ_N}{\goth W}_{k', k}^{m, n}, & {\goth W}_{k', k}^{m, n} \cap {\goth W}_{k', k}^{m', n'} \neq \emptyset ~ ~{\rm  iff} ~ ~ m \equiv m' ~ {\rm or}~ n \equiv n' ; \\
{\goth W}_{k', k}^{m, n} \cap {\goth W}_{k', k}^{m, n'}= {\rm Sing}^0({\goth W}_{k', k}), & {\goth W}_{k', k}^{m, n} \cap {\goth W}_{k', k}^{m', n}= {\rm Sing}^1({\goth W}_{k', k}), 
\elea(LWdeg)
which form a $G_N$-orbit, represented by $G_N {\goth W}_{k', k}^{0, 0}$; for instance, when $(k', k)= (-1, 0)$ with ${\goth L}_{1, 0}$ in (\req(Lijmn)), the isotropy subgroup of $G_N$ at ${\goth W}_{-1, 0}^{0, 0}$ is generated by $U{\sf u}_1^{-1}, U^2{\sf u}_2^{-1}$ with the order $2, 2N$ respectively, where $U{\sf u}_1^{-1}, U^2{\sf u}_2^{-1}$ are the involution and $\omega^\frac{1}{2}$-rotation of ${\goth W}_{-1, 0}^{0, 0}$ respectively.

(II). Three set of horizontal lines in (\req(lines)) are mutually disjoint. The superintegrable rapidity lines in Remark of Proposition \ref{prop:Ferm} are the horizontal lines defined by the first set of equations, i.e.  lines in  $\cup_{m, n \in \ZZ_N} \{ {\goth L}_{0, 0}^{m, n}, {\goth L}_{1, 1}^{m, n} \} ~ (= G_N {\goth L}_{0, 0}^{0, 0})$  in (\req(Lijmn)). Note that the isotropy subgroup of $G_N$ at ${\goth L}_{0, 0}^{0, 0}$ is generated by $U{\sf u}_2^{-1}, {\sf u}_2^2$ of order $2, N$ respectively, with $U{\sf u}_2^{-1}= {\rm id.}$ on ${\goth L}_{0, 0}^{0, 0}$, and  
${\goth L}_{0, 0}^{0, 0}/\langle {\sf u}_2^2 \rangle = \Lambda$.
$\Box$ \par \noindent \vspace{.1in}

\subsection{The family of hyperelliptic curves in chiral Potts model \label{subs.hypell}}
In this subsection, we describe the family of hyperelliptic curves with the symmetry group $\ZZ_2 \times D_N$ related to $G_N$ in (\req(GNDih)). First we consider the action of $H$ in (\req(Hrl)) on the fibration ${\goth W}$ in (\req(WFerm)). By (\req(uVMaut) ), the generators ${\sf V}_0,{\sf V}_1$ of $H$ are the automorphisms $M^{(0)},M^{(1)}$ of ${\goth W}_{k',k}$ in (\req(Sym)). 
Since $H$ acts freely on fibers ${\goth W}_{k',k}$ in (\req(Wk'k)), the function field of the quotient curve is generated by the $H$-invariant functions $t, \lambda$:  
\bea(ll)
[t, 1] = [ab, cd],  &  [\lambda , 1] =[d^N, c^N] \in \PZ^1 ,
\elea(tlambda)
satisfying the equation of a hyperelliptic curve of genus $N-1$ with $\ZZ_2 \times D_N$ symmetry:
\bea(llll)
W_{k',k} \simeq {\goth W}_{k',k}/H: & t^N = \frac{(1- k' \lambda)(1-k' \lambda^{-1})}{k^2}, & (t, \lambda) \in (\PZ^1)^2 , &  (k' \neq 0, \pm 1, \infty), 
\elea(Wk'kt)
where the $(\ZZ_2 \times D_N)$-structure is provided by (\req(GN0lij)) (\req(GNDih)) via the representation (\req(GtWf)):
\bea(lll)
\theta : (t, \lambda ) \mapsto (\omega t, \lambda), & \sigma: (t, \lambda ) \mapsto (t, \frac{1}{\lambda} ),  & \iota: (t, \lambda ) \mapsto (\frac{1}{t}, \frac{1-k'\lambda}{k'-\lambda}).
\elea(hyDih)
Indeed, the equation (\req(Wk'kt)) defines a surface over $\Lambda$ with the $\ZZ_2 \times D_N$ symmetry (\req(hyDih)): 
\bea(ll)
W =\{ (t, \lambda ; k', k) \in (\PZ^1)^2 \times \Lambda ~ | ~ k^2 t^N = (1- k' \lambda)(1-k' \lambda^{-1}) \} = \bigcup_{k', k \in \Lambda} W_{k',k} \longrightarrow \Lambda .
\elea(Whydef)
The fibers over the degenerated parameters are the following rational curves:
\bea(ll)
W_{0, 1}= W_{0,- 1} &=  \{ (t, \lambda ) \in (\PZ^1)^2 | t^N = 1 ~ {\rm or} ~ \lambda = 0^{\pm 1}  \} ; \\
W_{1,0} &=  \{ (t, \lambda ) \in (\PZ^1)^2 | t^N = \infty ~  {\rm or} ~ (\lambda - 1)^2=0  \}; \\
W_{- 1,0} &=  \{ (t, \lambda ) \in (\PZ^1)^2 | t^N = \infty ~  {\rm or} ~ (\lambda + 1)^2=0 \}; \\
W_{\infty,\infty_+ }= W_{\infty,\infty_- } &=  \{ (t, \lambda ) \in (\PZ^1)^2 | t^N = -1 ~ {\rm or} ~ \lambda = 0^{\pm 1}  \}.  
\elea(Whydeg)
Note that $\theta,  \sigma, \iota$ in (\req(hyDih)) preserve the degenerated curves in (\req(Whydeg)), inducing the automorphisms of $W_{0, \pm 1}, W_{\infty,\infty_\pm }$, but not for $W_{\pm 1,0}$. By the relation of $\kappa'$ and $k, k$ in (\req(kk'kap)), the variables $k'$, $k$, $\frac{1}{k}$ are the local coordinate of $\Lambda$ near $(k', k)= (0,\pm 1), (\pm 1,0), (\infty,\infty_{\pm})$ respectively. By which, $W$ in (\req(Whydef)) is a singular hypersurface in $(\PZ^1)^2 \times \Lambda$ with the singular locus 
\bea(ll)
{\rm Sing}(W)= \cup W_{\pm 1,0}, 
\elea(SingW)
where the singularity structure can be determined by the local coordinates $k, \epsilon ~ (=\lambda \mp 1), s ~(= t^{-1}$. In particular, $W$ is defined by the local equation, $k^2 = u  \epsilon^2$ near $\lambda = \pm 1, ~ t \neq 0$, or 
$k^2 = u  s^N $ near $t = \infty , ~ \lambda \neq \pm 1 $, where $u$ is a non-vanishing local function.

We now consider the action of $H$ on ${\goth W}$ in (\req(Wk'k)), and study the degeneration of curves in (\req(Wk'kt)) as $(k', k)$ tends to degenerated parameters. Denote the $H$-fixed point set of ${\goth W}$ by
$$
{\goth W}^H = \bigcup_{{\sf h} \in H, {\sf h} \neq {\rm id}} {\goth W}^{\sf h} , ~ ~ ~ {\goth W}^{\sf h} := \{ p \in {\goth W} ~ | ~ ~ {\sf h} ( p) = p \}.
$$
\begin{lem}\label{lem:Ffix} (i) ${\goth W}^H$ is a finite set consisting of singularities of degenerated fibers in (\req(SWdeg)) or (\req(SW01)):
 \bea(lll)
{\goth W}^H = {\goth W}^{{\sf V}_0} \cup {\goth W}^{{\sf V}_1}  \cup {\goth W}^{{\sf V}_0{\sf V}_1^{-1}} \cup {\goth W}^{{\sf V}_0{\sf V}_1}, &
{\goth W}^{{\sf V}_0} = \cup {\rm Sing}( {\goth W}_{\infty, \infty_{\pm}}), \\
{\goth W}^{{\sf V}_1} = \cup {\rm Sing}( {\goth W}_{0, \pm 1}), ~ ~ 
{\goth W}^{{\sf V}_0{\sf V}_1}= \cup {\rm Sing}^0( {\goth W}_{\pm 1, 0}), &
{\goth W}^{{\sf V}_0{\sf V}_1^{-1}}= \cup {\rm Sing}^1( {\goth W}_{\pm 1, 0}),
\elea(HFixF) 
where ${\sf V}_0 {\sf V}_1: [a,b, c, d] \mapsto [\omega a, \omega ^{-1} b, c,  d]$, ${\sf V}_0 {\sf V}_1 ^{-1}: [a,b, c, d] \mapsto [a, b, \omega c, \omega ^{-1} d]$. Furthermore, each singular set in (\req(HFixF)) is stable under the action of $H$, and $\overline{\rm Sing}_{k', k} (:={\rm Sing}({\goth W}_{k', k})/H)$ consists of two elements, $\overline{\rm Sing}^0_{k', k}$ and $\overline{\rm Sing}^1_{k', k}$, with the $a$-value $0, 1$ respectively.

(ii) The quotient of degenerated fibers in (\req(Wdeg)) by $H$ consists of lines given by 
\bea(lll)
{\goth W}_{0, \pm 1}/H &= \bigcup_{n \in \ZZ_N} \overline{\goth W}^{0, n}_{0, \pm 1}, ~ (\overline{\goth W}^{0, n}_{0, \pm 1}\simeq  {\goth W}^{0, n}_{0, \pm 1}/\langle {\sf V}_1 \rangle) , \\
{\goth W}_{\infty, \infty_{\pm}}/H &= \bigcup_{n \in \ZZ_N} \overline{\goth W}^{0, n}_{\infty, \infty_{\pm}}, ~ ( \overline{\goth W}^{0, n}_{\infty, \infty_{\pm}}\simeq  {\goth W}^{0, n}_{\infty, \infty_{\pm}}/\langle {\sf V}_0 \rangle), \\
{\goth W}_{\pm 1, 0}/H &= \left\{ \begin{array}{ll}
\overline{\goth W}^{0, 0}_{\pm 1, 0} ~ (\simeq {\goth W}^{0, 0}_{\pm 1, 0}) &{\rm if} ~ N ~ {\rm odd},  \\
\overline{\goth W}^{0, 0}_{\pm 1, 0}\cup \overline{\goth W}^{0, 1}_{\pm 1, 0}, (\overline{\goth W}^{0, j}_{\pm 1, 0} \simeq {\goth W}^{0, j}_{\pm 1, 0}/\langle {\sf V}_0^\frac{N}{2}{\sf V}_1^\frac{N}{2} \rangle ) &{\rm if} ~ N ~ {\rm even},
\end{array} \right.
\elea(Wdeg/H)
with the only intersection at $\overline{\goth W}^{0, n}_{k', k} \cap \overline{\goth W}^{0, n'}_{k', k} = \overline{\rm Sing}_{k', k}$ for $n \neq n'$.
\end{lem}
{\it Proof.} $(i)$ follows from the expression of automorphisms in $H$ and the Fermat relation of ${\goth W}$ in (\req(Wk'k)). The irreducible components of a degenerated fiber ${\goth W}_{k', k}$ are ${\goth W}_{k', k}^{m, n} ~ (m, n \in \ZZ_N)$ in (\req(LWdeg)), interchanged by $H$ via the relations:
$$
\begin{array}{llll}
{\sf V}_0: &{\goth W}_{0, \pm 1}^{m, n}, {\goth W}_{\infty, \infty_{\pm}}^{m, n}, {\goth W}_{\pm 1, 0}^{m, n} & \longrightarrow &{\goth W}_{0, \pm 1}^{m-1, n+1}, {\goth W}_{\infty, \infty_{\pm}}^{m, n}, {\goth W}_{\pm 1, 0}^{m-1, n-1}; \\
{\sf V}_1: &{\goth W}_{0, \pm 1}^{m, n}, {\goth W}_{\infty, \infty_{\pm}}^{m, n}, {\goth W}_{\pm 1, 0}^{m, n} & \longrightarrow &{\goth W}_{0, \pm 1}^{m, n}, {\goth W}_{\infty, \infty_{\pm}}^{m-1, n+1}, {\goth W}_{\pm 1, 0}^{m-1, n+1}.
\end{array}
$$
Hence ${\goth W}_{k', k}/H $ can be represented by the $H$-quotient classes in (\req(Wdeg/H)).
$\Box$ \par \noindent \vspace{.1in} 
All quotient curves, ${\goth W}_{k',,k}/H$ in (\req(Wk'kt)) and (\req(Wdeg/H)), form a family of curves over $\Lambda$: 
\bea(ll)
\pi: {\goth W}/H = \bigcup_{(k',k) \in \Lambda} {\goth W}_{k',k}/H \longrightarrow \Lambda .
\elea(HtGN)
By Lemma \ref{lem:Ffix} $(i)$, $H$ acts on ${\goth W}$ freely outside singularities of degenerated fibers, hence ${\goth W}/H$ is an orbifold with the singular locus
\bea(l)
{\rm Sing}({\goth W}/H) = \{ \overline{\rm Sing}_{0, \pm 1} , \overline{\rm Sing}_{\pm 1, 0} , \overline{\rm Sing}_{\infty, \infty_\pm} \}. 
\elea(SingWH)
\begin{lem}\label{lem:singWH} The local structure of ${\goth W}/H$ near a singularity $\overline{p}$ is given by 
\bea(ll)
({\goth W}/H , \overline{p}) \simeq \left\{\begin{array}{ll}
(\CZ^2/ \langle {\rm dia}[\omega, \omega] \rangle , \vec{0}), &{\rm for} ~ \overline{p}  \in \cup \overline{\rm Sing}_{0, \pm 1} \cup \overline{\rm Sing}_{\infty, \infty_\pm}, \\
(\CZ^2/ \langle {\rm dia}[\omega, \omega^{-1}] \rangle , \vec{0}), &{\rm for} ~\overline{p} \in \cup \overline{\rm Sing}_{\pm 1, 0}. 
\end{array} \right.
\elea(OrbWH)
\end{lem}
{\it Proof.} For $\overline{p} \in {\rm Sing}({\goth W}/H)$, let $p$ be an element in ${\rm Sing}( {\goth W}_{k', k})$ with $\overline{p} = H \cdot p$. Then $p$ is fixed by an automorphism ${\sf W}$ in (\req(HFixF)) with ${\sf W}= {\sf V}_0, {\sf V}_1$ or ${\sf V}_0{\sf V}_1^{\pm}$.  The local structure of ${\goth W}/H$ near $\overline{p}$ is isomorphic to ${\goth W}/\langle {\sf W} \rangle$ near $p$. By (\req(SWdeg)), among the four coordinates $a, b, c, d$, two are zeros at $p$, which provide the local coordinate system near $p$ in ${\goth W}$. By (\req(HFixF)), the expression of ${\sf W}$ gives rise to the local structure of $({\goth W}/\langle {\sf W} \rangle, p)$, then follows (\req(OrbWH)).
$\Box$ \par \noindent \vspace{.1in} 
Note that the relation (\req(tlambda)) defines a birational correspondence between ${\goth W}/H$ in (\req(HtGN)) and $W$ in (\req(Whydef)) over $\Lambda$:
\bea(lll)
\varrho: {\goth W}/H \rightharpoonup  W , & H \cdot q \mapsto (t, \lambda; k' k) =(\frac{ab}{cd}, \frac{d^N}{c^N}; \pi (q)) &   ( q =[a,b, c, d] \in {\goth W}),
\elea(birWHW)
which allows us to make the identification (\req(Wk'kt)) for smooth fibers, with the birational or two-one equivalences between $H$-classes in (\req(Wdeg/H)) and  components of the degenerated fibers (\req(Whydeg)):
\bea(ll)
{\goth W}_{0, \pm 1}/H \stackrel{\sim}{\rightharpoonup} (W_{\pm 1, 0})_{t^N = 1}, ~ {\goth W}_{\infty, \infty_{\pm}}/H \stackrel{\sim}{\rightharpoonup} (W_{\infty, \infty_{\pm}})_{t^N = -1},  &{\goth W}_{\pm 1, 0}/H \stackrel{2:1}{\rightharpoonup}  (W_{\pm 1, 0})_{\lambda = \pm 1}. 
\elea(degcor)
By (\req(SWdeg)), the fundamental locus of the birational correspondence (\req(birWHW)), where $\varrho$ is not well-defined,  consists of 10 singularities of ${\goth W}/H $:
$$
ab=cd=0 : \overline{\rm Sing}_{0, \pm 1} , \overline{\rm Sing}_{\infty, \infty_\pm} , ~  c=d = 0 : \overline{\rm Sing}^1_{\pm 1, 0}.
$$
In order to replace $\varrho$ by a regular morphism, we consider the minimal resolution $\widehat{{\goth W}_H}$ of ${\goth W}/H$, 
\bea(ll)
 \widehat{\eta} : \widehat{{\goth W}_H}  \longrightarrow {\goth W}/H , & \widehat{\varrho}~ (= \varrho \cdot \widehat{\eta}):  \widehat{{\goth W}_H} \longrightarrow W .
\elea(resl)
It is known in algebraic geometry that $\widehat{\eta}$ is a regular morphism from the non-singular variety $\widehat{{\goth W}_H}$ to ${\goth W}/H$ where ${\rm Sing}({\goth W}/H)$ in  ${\goth W}/H$ are replaced by exceptional curves $\widehat{\eta}^{-1}({\rm Sing}({\goth W}/H))$ in $\widehat{{\goth W}_H}$. The relationship between $\widehat{{\goth W}_H}$ and the resolution of $W$ is described by the following proposition.
\begin{prop}\label{prop:relW} (i) For even $N$, $\widehat{\varrho}$ is a regular birational morphism, by which $\widehat{{\goth W}_H} $ is a resolution of $W$ via $\widehat{\varrho}$, i.e. $\widehat{\varrho}$ is a regular birational morphism, biregular between $\widehat{{\goth W}_H} \setminus \widehat{\varrho}^{-1}({\rm Sing}(W)) $ and $W \setminus {\rm Sing}(W)$.

(ii) For odd $N$, the fundamental locus of $\widehat{\varrho}$ consists of two elements $\widehat{o}_{\pm 1, 0}$ lying over $\overline{\rm Sing}^1_{\pm 1, 0}$ respectively. The blow-up of $\widehat{{\goth W}_H}$ centered at $\widehat{o}_{\pm 1, 0}$ is a resolution of $W$ via the $\widehat{\varrho}$-induced regular birational morphism. 
\end{prop}
{\it Proof.} By (\req(SingW)), ${\rm Sing}(W)= \cup W_{\pm 1,0}$, with normal-crossing double curves outside $t =\infty, 0$ in (\req(Whydeg)). The normalization of $W$ outside $t= \infty, 0$ in  $W_{\pm 1,0}$ provides the resolution, which can be identified with ${\goth W}/H$ outside ${\rm Sing}({\goth W}/H)$ via $\varrho$ in (\req(birWHW)). Hence we need only to study the relationship between $W$ and the minimal resolution of  ${\goth W}/H$ near ${\rm Sing}({\goth W}/H)$ locally. Since ${\goth W}/H$ possesses only  orbifold-singularities (\req(OrbWH)), its minimal resolution is constructed by either Hirzebruch-Jung continued fraction method \cite{Hir, Jun} or techniques in toric geometry \cite{Oda}. Consider the first type of singularity $\overline{p}$ in (\req(OrbWH)) with the local coordinates $(z_1, z_2) ~ (\in \CZ^2) $ of ${\goth W}$ near $p \in {\rm Sing}({\goth W}_{k', k})$, where $(z_1, z_2) = (\frac{a}{c}, \frac{d}{b}), (\frac{b}{d}, \frac{c}{a}),(\frac{a}{d}, \frac{c}{b}), (\frac{b}{c}, \frac{d}{a}) $ for $p \in \cup {\rm Sing}^0({\goth W}_{0, \pm 1})$, $\cup {\rm Sing}^1({\goth W}_{0, \pm 1})$, $\cup {\rm Sing}^0({\goth W}_{\infty, \infty_\pm})$, $\cup {\rm Sing}^1({\goth W}_{\infty, \infty_\pm})$ respectively. Then the exceptional curve $E_{\overline{p}} ~ (= \widehat{\eta}^{-1} (\overline{p}))$  in the minimal resolution $\widehat{{\goth W}_H}$ is a rational curve $E_{\overline{p}}\simeq \PZ^1$  with self-intersection number equal to $-N$. Indeed, the smooth manifold $\widehat{{\goth W}_H}$ near $E_{\overline{p}}$ is covered by two charts with local coordinates $(z_1^N, \frac{z_2}{z_1})$ or $(\frac{z_1}{z_2}, z_2^N)$, and the projective coordinate of  $E_{\overline{p}}$ can be identified with the $t$-variable in (\req(birWHW)): $[z_1: z_2] = [ab, cd]$. Near $E_{\overline{p}}$ in $\widehat{{\goth W}_H}$, the local defining equation of $\widehat{\eta}^{-1}({\goth W}_{k', k})$ is given by
$$
\begin{array}{llll}
(\frac{z_1}{z_2})^N = t^N =  1, & {\rm or} & z_1^N = z_2^N = \pm \lambda = 0 ,  & p \in \cup {\rm Sing}^0({\goth W}_{0, \pm 1}) ; \\
(\frac{z_1}{z_2})^N = t^N =  1, & {\rm or} & z_1^N = z_2^N = \pm \lambda^{-1}=0, & p \in \cup {\rm Sing}^1({\goth W}_{0, \pm 1}) ; \\
(\frac{z_1}{z_2})^N = t^N = -1 , & {\rm or} & z_1^N = - z_2^N = \mp {\rm i} \lambda^{-1}=0 , &  p \in \cup {\rm Sing}^0({\goth W}_{\infty, \infty_\pm}) ; \\
(\frac{z_1}{z_2})^N = t^N = -1, & {\rm or} & z_1^N = - z_2^N = \mp {\rm i} \lambda =0, &  p \in \cup {\rm Sing}^1({\goth W}_{\infty, \infty_\pm}) ,
\end{array}
$$
which are the same as the corresponding $W_{k', k}$ in (\req(Whydeg)). By (\req(SingW)), $\widehat{\varrho}$ in (\req(resl)) defines a biregular isomorphism between $ \widehat{{\goth W}_H} \setminus \widehat{\eta}^{-1}({\goth W}_{\pm 1, 0})$ and $W \setminus {\rm Sing}(W)$. Hence $\widehat{\varrho}$ is a resolution of $W$ outside ${\rm Sing}(W)$ in (\req(SingW))\footnote{The local charts near $E_{\overline{p}}$ in $\widehat{{\goth W}_H}$ is biregular equivalent to $\lambda^{\pm 1}=0$ in $W$ locally only, but not the entire lines $t^N = \pm 1$.  For example, in ${\rm Sing}^0({\goth W}_{0, \pm 1})$ case, the local biregular morphism is given by $z_1^N = \frac{\lambda - k'}{k}, z_2^N= \frac{k \lambda}{1- k\lambda}$, and  the equation (\req(Whydef)) of $W$ with $t^N =1$ becomes $k'(\lambda - (k' +{\rm i} k))( \lambda - (k' -{\rm i} k))=0$. Hence $t^N = \pm 1$ are $N$-lines in $W$, each with self-intersection number $-2$.}. It remains to consider the local behavior of  $\widehat{\varrho}$ near exceptional curves over the second type of singularity $\overline{p}$ in (\req(OrbWH)), where the local coordinate $(z_1, z_2)$ of ${\goth W}$ near $p \in {\rm Sing}({\goth W}_{\pm 1, 0})$ is given by
$(z_1, z_2) = (\frac{a}{d}, \frac{b}{c}), (\frac{c}{b}, \frac{d}{a})$ for $p \in \cup {\rm Sing}^0({\goth W}_{\pm 1, 0})$, $\cup {\rm Sing}^1({\goth W}_{\pm 1, 0})$ respectively.  The orbifold singularity at $\overline{p}$ is of type $A_{N-1}(= \langle {\rm dia}[\omega, \omega^{-1}] \rangle)$ . It is known that the minimal resolution of $\CZ^2/A_{N-1}$ is covered by $N$ charts 
with the coordinates system provided by toric geometry (for example, as illustrated in \cite{CR} section 3): 
\bea(ll)
{\cal U}_j \simeq \CZ^2 : (u_j, v_j) = (z_1^{j+1}z_2^{1-N+j},z_1^{-j} z_2^{N-j}), ~ \widehat{o}_j : (u_j, v_j)=(0, 0) ,
\elea(ResA)
for $j=0, \ldots, N-1$. The exceptional divisor is $E_{\overline{p}} = E_1+ \cdots + E_{N-1}$, where $E_j$ is the rational $(-2)$-curve in the minimal resolution joining $\widehat{o}_{j-1}$ to $\widehat{o}_j$ defined by $v_{j-1}= 0 = u_j$. The divisors $D_0, D_N$  defined by $u_0=0, v_{N-1}=0$ in ${\cal U}_0$ or ${\cal U}_{N-1}$ are the proper transform of $(z_1^N=0)$ or $(z_2^N=0)$ in $\CZ^2/A_{N-1}$ respectively. First, we consider the singularity $\overline{p} = \overline{\rm Sing}^0_{\pm 1, 0}$, where  $\varrho$ in (\req(birWHW)) is well-defined. In the orbifold model (\req(OrbWH)), ${\goth W}_{\pm 1, 0}$ near $p$ corresponds to $N$-lines in $\CZ^2$: $z_1^N= - z_2^N$, whose quotient curve in $\CZ^2/ A_{N-1}$ corresponds to ${\goth W}_{\pm 0}/H$ near $\overline{p}$ in  (\req(Wdeg/H)). 
The inverse process of normalization of $W$ near $(t, \lambda)=(0, \pm 1) \in W_{\pm 1, 0}$ in (\req(Whydeg)) is equivalent to the identification of $N$-lines $z_1^N= - z_2^N$ in $\CZ^2$ via the following local automorphisms:
$$
(z_1, z_2) \mapsto \left\{ \begin{array}{ll}
(\omega z_1, \omega^{-1} z_1), (-z_1, -z_2), & {\rm for ~ odd ~} N, \\
(\omega z_1, \omega^{-1} z_1), (-z_2, -z_1), & {\rm for ~ even ~} N .
\end{array}
\right.
$$ 
Note that the above automorphisms are the restriction of ${\sf V}_0{\sf V}_1, T^{2N}, U{\bf u}_2^{N-1} \in \widetilde{G}_N $ on ${\goth W}_{\pm 1, 0}$ respectively, with $U{\bf u}_2^{N-1} = {\sf V}_0^\frac{1-N}{2}{\sf V}_1^\frac{1-N}{2} $ for odd $N$ and $T^{2N}= {\sf V}_0^\frac{N}{2}{\sf V}_1^\frac{N}{2}$ for even $N$. The induced identification of curves in the orbifold corresponds to the $(2:1)$-morphism in (\req(degcor)). Therefore the minimal resolution over $\overline{\rm Sing}^0_{\pm 1, 0}$ provides a resolution of $W$ near $(t, \lambda)=(0, \pm 1) \in W_{\pm 1, 0}$.  
Next, we examine the behavior of $\widehat{\varrho}$ in (\req(resl)) near the exceptional divisor $E_{\overline{p}}$ in the minimal resolution $\widehat{\CZ^2/A}_{N-1}$ of $\CZ^2/A_{N-1}$ for $\overline{p} \in \overline{\rm Sing}^1_{\pm 1, 0}$ . Since the local coordinates $(z_1, z_2)$ of ${\goth W}$ near $p$ are related to the coordinates of $t, \lambda (= \frac{d^N}{c^N})$ by $t^{-1} = z_1z_2 $, $(\lambda -k')z_2^N = \lambda  (1 -k'\lambda) z_1^N$, the rational map $\varrho$ in (\req(birWHW)) at $k'= \pm 1$ is defined by  $t = \infty , ~ [\mp \lambda, 1] = [z_2^N, z_1^N]$. In the affine chart ${\cal U}_j $ in (\req(ResA)), $z_1^N= u_j^{N-j} v_j^{N-1-j}, z_2^N= u_j^j v_j^{j+1}$, by which $z_1^N, z_2^N$ give rise to two regular functions of $\widehat{\CZ^2/A}_{N-1}$ with the zero-divisor $
(z_1^N=0) = \sum_{j=1}^{N-1} (N-j) E_j + N D_0$, $(z_2^N=0)= \sum_{j=1}^{N-1} j E_j + N D_N$,
whose ratio  defines a rational map 
\bea(ll)
\chi: \widehat{\CZ^2/A}_{N-1} \longrightarrow \PZ^1 ,& * \mapsto [z_1^N, z_2^N](*).
\elea(chi)
Then $\chi (D_0)= [0, 1], \chi (D_N)= [1, 0]$, and  $\chi (E_j)= [0, 1]$ if $j < \frac{N}{2}$, and $\chi (E_j)= [1, 0]$ if $ j > \frac{N}{2}$.  For even $N$, $\chi$ in (\req(chi)) is indeed a regular birational morphism, which induces a double cover  of $E_\frac{N}{2}$ over $\PZ^1$ branched at $0, \infty$. Then follows $(i)$ with $E_\frac{N}{2}$ as a double cover of the $\PZ^1$-curve in (\req(Whydeg)) defined by $t^N= \infty$.
When $N$ is odd, the fundamental locus  of $\chi$ in (\req(chi)) consists of only one element: $\widehat{o}_\frac{N-1}{2}$. Indeed, one finds $[z_1^N, z_2^N]=[u_\frac{N-1}{2},   v_\frac{N-1}{2}]$ in the affine chart ${\cal U}_\frac{N-1}{2}$ centered at $\widehat{o}_\frac{N-1}{2}$. Then $\chi$ can be lifted to a regular morphism from the blow-up of $\widehat{\CZ^2/A}_{N-1}$ at $\widehat{o}_\frac{N-1}{2}$ to $\PZ^1$, which is identified with the exceptional curve $E$ in the blow-up manifold. Since $\widehat{\varrho}$ in (\req(resl)) is equivalent to $\chi$ in (\req(chi)) near the exceptional divisor $E_{\overline{p}}$,  $(ii)$ follows with $\widehat{o}_\frac{N-1}{2}= \widehat{o}_{\pm 1, 0}$, and $E$ isomorphic to the $(t^N= \infty)$-curve in (\req(Whydeg)).
$\Box$ \par \vspace{.1in} 

As in (\req(Whydef)),(\req(HtGN)) and (\req(birWHW)), we also consider the quotients of ${\goth W}$ by $H_r, H_l$ in (\req(Hrl)):
\bea(ll)
\pi_r: {\goth W}/H_r  \longrightarrow \Lambda , & \pi_l: {\goth W}/H_l \longrightarrow \Lambda ,
\elea(HrlGN)
which are related to the following families of hyperelliptic curves with $\ZZ_2 \times D_N$ symmetry over $\Lambda$:
\bea(ll)
W_r =\{ (t_r, \lambda_r ; k', k) \in (\PZ^1)^2 \times \Lambda ~ | ~ t_r^N = (k- {\rm i} k' \lambda_r)(k-{\rm i} k' \lambda_r^{-1}) \} \longrightarrow \Lambda , \\
 W_l =\{ (t_l, \lambda_l ; k', k) \in (\PZ^1)^2 \times \Lambda ~ | ~ k^{' 2} t_l^N = (1- k \lambda_l )(1-k \lambda_l^{-1}) \} \longrightarrow \Lambda . 
\elea(Wrl)
via the birational correspondence (see, \cite{R04} (27))\footnote{The variables $(T_r, \Lambda_r), (T_l, \Lambda_l)$ in \cite{R04} (27) are related to $(t_r, \lambda_r), (t_l, \lambda_l)$  here by $(T_r, \Lambda_r)= (t_r, \lambda_r)$, $(T_l, \Lambda_l)= (t_l^{-1}, \frac{1-k \lambda_l}{k- \lambda_l})$. Indeed, the relation between $(T_l, \Lambda_l)$ and $(t_l, \lambda_l)$ corresponds to the $W_{k, k'}$-automorphism $\iota$ in (\req(hyDih)), which is induced by the automorphism $U{\sf u}_2^{-1}{\sf u}_1$ of ${\goth W}_{k, k'}$.}:
\bea(ll)
\varrho_r: {\goth W}/H_r \rightharpoonup  W_r , & H_r \cdot q \mapsto (t_r, \lambda_r; k', k) =(\frac{ac}{bd}, \frac{{\rm i}d^N}{b^N}; \pi_r (q));   \\
\varrho_l: {\goth W}/H_l \rightharpoonup  W_l , & H_l \cdot q \mapsto (t_l, \lambda_l; k', k) =(\omega^{-1/2} \frac{ad}{bc}, \frac{b^N}{c^N}; \pi_l (q)), \\
\elea(birWrl)
where $q =[a,b, c, d] \in {\goth W}$. Note that for $N \geq 3$, the above fibration are different from (\req(birWHW)) by Lemma \ref{lem:HsST}. However, since $(STS^{-1}) H_r = H (STS^{-1}), S H_l = HS$ by (\req(GNlr)), one finds the isomorphic relations between (\req(birWHW)) and (\req(birWrl)), using (\req(STWk')) and the representation of $STS^{-1}, S$ in (\req(GtWf)): 
\bea(clll)
STS^{-1}&: {\goth W}/H_r \simeq {\goth W}/H, & {\goth W}_{k', k}/H_r \simeq {\goth W}_{\frac{{\rm i} k'}{k}, \frac{1}{k}}/H, & (H_r \cdot q \mapsto  H \cdot STS^{-1} q ), \\
&\varrho_r = \varrho (STS^{-1}),&  W_{r; k', k} \simeq W_{\frac{{\rm i} k'}{k}, \frac{1}{k}}, &(t_r, \lambda_r) \stackrel{STS^{-1}}{\mapsto ~ ~ } (t, \lambda)= (t_r, \lambda_r);  \\
S &: {\goth W}/H_l \simeq {\goth W}/H, & {\goth W}_{k', k}/H_l \simeq {\goth W}_{k, k'}/H, &(H_l \cdot q \mapsto  H \cdot Sq) , \\
 & \varrho_l= \varrho S & W_{l; k', k} \simeq W_{k, k'}, &(t_l, \lambda_l) \stackrel{S}{\mapsto} (t, \lambda)= (t_l, \lambda_l).
\elea(isHrlH)
Through the equivalences in (\req(isHrlH)), one can derive the relationship between the minimal resolution of ${\goth W}/H_r$ and $W_r$, or ${\goth W}/H_l$ and $W_l$ in (\req(birWrl)), from Lemma \ref{lem:singWH} and Proposition \ref{prop:relW}.

On the other hand, the normalizer of $H$ in (\req(NH)) gives rise to an automorphism group of ${\goth W}/H$:
\bea(ll)
N(H)/H = \langle G_N, T^{-1}, ST^2 S^{-1} \rangle/H =  \langle \theta, \sigma, \iota , T^{-1}, ST^2 S^{-1} \rangle , & T^4 = \theta^2, ~ ~ ST^4 S^{-1}= 1,
\elea(autWH)
which induces an action on $\widehat{\goth W}_H$ and $W$ in (\req(resl)). Using (\req(GtWf)), one finds the expression of the $N(H)/H$-action on $W$. In particular, the $\langle T^{-1}, ST^2 S^{-1} \rangle$-expression of $W$ is given by
\bea(lll)
T:& W_{k', k} \longrightarrow W_{\frac{1}{k'}, \frac{{\rm i}k}{k'} }, & (t, \lambda; k',k) \mapsto 
(\omega^\frac{1}{2} t, \lambda^{-1}; \frac{1}{k'}, \frac{{\rm i}k}{k'}), \\
T^2:& W_{k', k} \longrightarrow W_{k', -k }, & (t, \lambda; k',k) \mapsto 
(\omega t, \lambda ; k', -k ), \\
ST^{-2} S^{-1}:& W_{k', k} \longrightarrow W_{-k', k}, &(t, \lambda; k',k) \mapsto (t, - \lambda; -k', k),
\elea(TST2SW)
with the induced  relationship of $(\ZZ_2 \times D_N)$-symmetry inherited from those of $\widehat{\goth W}$ in (\req(NHUij)); for example, $T^{-1}(\theta, \sigma, \iota) T = (\theta, \sigma, \theta^{-1} \sigma \iota)$ corresponds to the $T^{-1}$-relation in (\req(NHUij)). Note that $T: {\goth W}_{k', k} \simeq {\goth W}_{\frac{1}{k'}, \frac{{\rm i}k}{k'}}$ in (\req(GtWf)), and the automorphism $T$ in (\req(TST2SW)) is the rapidity-identification in the (Kramers-Wannier) duality of chiral Potts model (\cite{KW}, \cite{R905} (3.9) (3.12)).

\subsection{Fermat K3 surface and elliptic fibration \label{subs.K3}}
When $N=2$, the symmetries discussed in Subsections \ref{subs.FSym} and \ref{subs.hypell} are indeed  the elliptic and modular symmetries of elliptic curves expressed by the theta functions of half-integer characteristics:
$$
\begin{array}{lllll}
 \vartheta_1  ( v, \tau ) &(= \vartheta^{(\frac{-1}{2}, \frac{-1}{2})} ( v , \tau )) & = 
2 q_0 q^{ \frac{1}{8}} \sin \pi v \prod_{n=1}^{ \infty} (
1 - 2 q^n \cos 2 \pi v + q^{2n} ) \\
 \vartheta_2 ( v  , \tau )  & (=  
\vartheta^{(\frac{1}{2}, 0)} ( v , \tau ) )& = 
2 q_0 q^{ \frac{1}{8}} \cos \pi v \prod_{n=1}^{ \infty} (
1 + 2 q^n \cos 2 \pi v + q^{2n} ) \\
 \vartheta_3 ( v  , \tau )  &(=  
\vartheta^{(0, 0)} ( v , \tau ) )& =  
q_0 \prod_{n=1}^{ \infty} ( 1 + 2 q^{n-\frac{1}{2}} \cos 2 \pi v + 
q^{2n-1} ) \\
 \vartheta_4 ( v, \tau )  &(=  
\vartheta^{(0, \frac{1}{2})} (v, \tau ) )& =  q_0  \prod_{n=1}^{ \infty} (
1 - 2 q^{n - 1/2} \cos 2 \pi v + q^{2n - 1} ), \\
\end{array}
$$
where $ q= e^{2 \pi {\rm i} \tau} , q_0 : = \prod_{n=1}^{\infty} ( 1 - q^n )$, $\tau \in \HZ$ (the upper-half plane). The above theta functions  satisfy the elliptic and modular properties:
\bea(ll)
\vartheta_{1, 2} (v +1 , \tau) = - \vartheta_{1, 2} (v, \tau), & \vartheta_{1, 4} ( v + \tau, \tau ) = - e^{- \pi {\rm i} (\tau + 2 v)} \vartheta_{1, 4} ( v , \tau );  \\
\vartheta_{3, 4} (v +1 , \tau) =  \vartheta_{3, 4} (v, \tau), & \vartheta_{2, 3} ( v + \tau , \tau ) = e^{- \pi {\rm i} (\tau + 2 v)} \vartheta_{2, 3} ( v , \tau );  \\
\vartheta_{1, 2}  ( v, \tau +1 ) = e^\frac{\pi {\rm i}}{4}\vartheta_{1, 2}  ( v, \tau ), & \vartheta_{3, 4}  ( v, \tau +1 ) = \vartheta_{4, 3}  ( v, \tau ), \\
\vartheta_1  ( \frac{v}{\tau}, \frac{-1}{\tau}) = -{\rm i} (-{\rm i} \tau)^\frac{1}{2} e^\frac{{\rm i} \pi v^2}{\tau} \vartheta_1  ( v, \tau), & \vartheta_{2, 3, 4}  ( \frac{v}{\tau}, \frac{-1}{\tau}) = (-{\rm i} \tau)^\frac{1}{2} e^\frac{{\rm i} \pi v^2}{\tau} \vartheta_{4, 3, 2}  ( v, \tau).
\elea(ellmod/2) 
and the algebraic relations:
\bea(ll)
 \vartheta_1(v, \tau)^2   = k \vartheta_4(v, \tau)^2 -k' \vartheta_2 (v, \tau)^2, & \vartheta_2(v, \tau)^2   = k \vartheta_3(v, \tau)^2  - k' \vartheta_1 (v, \tau)^2 ,  \\
\vartheta_3(v, \tau)^2 = k \vartheta_2(v, \tau)^2 + k' \vartheta_4 (v, \tau)^2 ,&
\vartheta_4(v, \tau)^2 = k\vartheta_1(v, \tau)^2 +k' \vartheta_3 (v, \tau)^2    ,
\elea(ell2alg)
where $k= \frac{\vartheta_2(0)^2}{\vartheta_3(0)^2} , k' = \frac{\vartheta_4(0)^2}{\vartheta_3(0)^2} $ with $k^2 + k'^2 = 1$. The first two relations in (\req(ellmod/2)) yield $\vartheta_j (v +2 , \tau)= \vartheta_j (v, \tau)= 
e^{4\pi {\rm i} (\tau + v)} \vartheta_j (v +2\tau , \tau)$, hence by (\req(ell2alg)), we find the uniformization of the elliptic curve ${\goth W}^{(2)}_{k', k}$\footnote{Here the parameterization differs from those in \cite{B93c} section 3 or \cite{R905} section 3.4 by some minus signs, where $[a, b, c, d]= [-\vartheta_1(v ; \tau ), - \vartheta_2(v ; \tau ), \vartheta_3(v ; \tau ), \vartheta_4 (v ; \tau )]$ with $v= \frac{u}{2I}, I = \frac{\pi}{2} \vartheta_3^2(0, \tau )$.}:
\bea(llll)
\CZ/(2\ZZ+ 2\tau \ZZ) \simeq {\goth W}^{(2)}_{k', k}, & [v] \mapsto & [a, b, c, d]:= [\vartheta_1, \vartheta_2, \vartheta_3, \vartheta_4](v, \tau), & k' \neq 0, \pm 1, \infty.
\elea(ellCP2)
The family ${\goth W}^{(2)}$ in (\req(WFerm)) for $N=2$ is a elliptic K3 surface over $\Lambda ~(\simeq \HZ/PSL_2(\ZZ_4))$. The generators of $\widetilde{G}_2$ for the fibration in (\req(GtWf)) can be identified with 
\bea(lll)
{\sf u}_1:& [a, b, c, d] \mapsto  &[{\rm i} d, c,  b, {\rm i} a] = [a, b, c, d] (v +\frac{\tau}{2}, \tau), \\
{\sf u}_2: &[a, b, c, d] \mapsto & [b, -a,  d, c] = [a, b, c, d] (v +\frac{1}{2}, \tau) , \\
U : & [a, b, c, d] \mapsto &[-a, b,  c, d] =  [a, b, c, d] ( -v , \tau) , \\
S: &  [a, b, c, d] \mapsto & [-{\rm i} a, d, c, b]= [a, b, c, d] ( \frac{v}{\tau}, \frac{-1}{\tau}), \\
T:  &  [a, b, c, d] \mapsto & = [{\rm i}^\frac{1}{2} a, {\rm i}^\frac{1}{2} b, d, c]= [a, b, c, d](v, \tau+1 ) ,
\elea(G2til)
in which $\langle {\sf u}_1, {\sf u}_2, U \rangle$ is equal to the automorphism group $G_2$ of (\req(ellCP2)). Indeed, the K3 surface ${\goth W}^{(2)}$ is parametrized by $\CZ \times \HZ$, with the action of
the semi-product group $\RZ^2 \ast SL_2(\RZ)$ from the right, where the conjugation of ${\rm M}=  \left( \begin{array}{cc}
A&B\\ 
C&D
\end{array}
\right) \in SL_2(\RZ)$ on $u = ({}^a_b) \in \RZ^2$ in $\RZ^2 \ast SL_2(\RZ)$ is given by ${\rm M} \cdot u \cdot {\rm M}^{-1} = {\rm M}^t u$ (matrix multiplication). It is known that the $\RZ^2 \ast SL_2(\RZ)$-action gives rise to a left-action on entire functions $\varphi ( v, \tau)$  of $\CZ \times \HZ$: $ ({\rm g}|\varphi) ( v, \tau)= \lambda (v, \tau ; {\rm g})\varphi \bigg( ( v, \tau){\rm g}\bigg)$ where $\lambda (v, \tau ; u)= e^{ \pi {\rm i} (a^2 \tau + 2 a( v + b))}$ and $\lambda (v, \tau ; {\rm M}) = (C\tau+D)^{\frac{-1}{2}}
e^{ \frac{-\pi {\rm i}C v^2 }{C \tau+D}}$\footnote{In this paper, the automorphism group acts on spaces from the right with the induced left-action on functions, different from the convention in \cite{R99} where the automorphism group acts on spaces from the left with induced operators on functions acting from the right. Indeed, the action of the semi-product $\RZ^2*SL_2(\RZ)$ on $\CZ\times \HZ$ from the left was given by formulas in \cite{R99} page 3067 with the conjugation of $M \in SL_2(\RZ)$ on $v \in \RZ^2$ given  by  $M^{-1} \cdot v \cdot M = vM$ (the matrix product). Changing the left action 
in \cite{R99} to the right action, one obtains the convention used in this paper: $M \cdot v^t \cdot M^{-1} = M^t v^t$. Similarly the left action $M |*$ on functions here is the same as the right action $* |M$ in \cite{R99}.}.  The $\widetilde{G}_2$-action of ${\goth W}^{(2)}$ in (\req(G2til)) is indeed induced from an action of $\CZ \times \HZ$ so that each $g \in \widetilde{G}_2$ can be identified with some ${\rm g} \in \QZ^2 \ast SL_2(\ZZ)$ acting on  $\CZ \times \HZ$. The $\widetilde{G}_2$-generators in (\req(G2til)) correspond to the following elements in $\QZ^2 \ast SL_2(\ZZ)$:
$$
\begin{array}{lll}
g= {\sf u}_1, {\sf u}_2,  S , T , U & \leftrightarrow & {\rm g}= u_1 = ({}^\frac{1}{2}_{0}) , u_2= ({}^0_{\frac{1}{2}}) ,  {\rm S}= \left( \begin{array}{cc}
0&-1\\ 
1&0
\end{array}
\right), {\rm T}= \left( \begin{array}{cc}
1&1\\ 
0&1
\end{array}
\right), {\rm S}^2 ; \\
S {\sf u}_i S^{-1}, T {\sf u}_i T^{-1}  & \leftrightarrow & {\rm S}^{-1} u_i {\rm S}, {\rm T}^{-1} u_i {\rm T} . \\
\end{array}
$$
Note that in the above correspondence,  $g_1 g_2 \in \widetilde{G}_2$ corresponds to ${\rm g}_2 {\rm g}_1 \in \QZ^2 \ast SL_2(\ZZ)$ for $g_i \in \widetilde{G}_2$ induced from ${\rm g}_i \in \QZ^2 \ast SL_2(\ZZ) ~ (i=1, 2)$. Hence we obtain the description of ${\goth W}^{(2)}$, $\widetilde{G}_2$ and $G_2$ in terms of uniformization of elliptic curves:
$$
\begin{array}{lll}
{\goth W}^{(2)} &\simeq (\CZ \times \HZ)/\langle \langle4u_1, 4u_2,  (2u_2) \ast {\rm T}^{-4}  \rangle \rangle_{\rm N} , \\
\widetilde{G}_2 &\simeq  (\ZZ_4^2 \ast SL_2(\ZZ)) /\langle \langle  ({}^0_2 ) \ast {\rm T}^{-4}  \rangle \rangle_{\rm N}, &
G_2 \simeq \ZZ_4^2 \ast \langle S^2 \rangle .
\end{array}
$$
{\bf Remark.} The above $G_2$ is a representation of relations in (\req(G2)). Even though there are infinity many symmetries for a single fiber (\req(ellCP2)), only those in $G_2$ can be extended to automorphisms of the fibration ${\goth W}^{(2)}$, including degenerated fibers in (\req(Wdeg)). By Proposition \ref{prop:lineF} Remark (II), the superintegrable rapidities in (\req(ellCP2)) is the $G_2$-orbit of $\frac{1}{4}$, consisting of 16 elements. 
$\Box$ \par \vspace{.1in}

We now consider the hyperelliptic family in Subsection \ref{subs.hypell} for $N=2$. By Lemma \ref{lem:HsST} $(i)$, the subgroups in (\req(Hrl)) are all equal, $H= H_r =H_l =G'_{2, 1} = \langle 2u_1, 2u_2 \rangle$ with $\widetilde{G}_2/G'_{2, 1} $ in (\req(G/G12)); so are the fibrations in (\req(HtGN)) and (\req(HrlGN)): ${\goth W}/H= {\goth W}/H_r = {\goth W}/H_l$ with only orbifold singularities at (\req(SingWH)), all of type $A_1$ by Lemma \ref{lem:singWH}. Hence the minimal resolution $\widehat{{\goth W}_H}$ of ${\goth W}^{(2)}/H $ in (\req(resl)) is a K3 surface. Indeed, $\widehat{{\goth W}_H}$ is also an elliptic fibration over $\Lambda$ with the fibers described by
$$
(\widehat{{\goth W}_H})_{k',k} = \left\{\begin{array}{ll} {\goth W}^{(2)}_{k', k}/H \simeq \CZ/\ZZ + \tau \ZZ,  &{\rm if} ~ k' \neq 0, \pm 1, \infty; \\
F_+ + F_- + E_+ + E_-, & {\rm if} ~ (k', k)=(0, \pm 1), (\pm 1, 0), (\infty, \infty_\pm), 
\end{array} \right.
$$
where the parameter $\tau$ is inherited from (\req(ellCP2)), and $F_\pm$ are the proper transforms of degenerated fibers in (\req(Wdeg/H)), $E_\pm$ the exceptional divisors, of which all are rational $(-2)$-curves with the intersection only between $E_\pm$ and $F_\pm$. By Proposition \ref{prop:relW} $(i)$, $\widehat{{\goth W}_H}$ can be regarded as a resolution of $W$ in (\req(Whydef)), so the same  for $W_r, W_l$ in (\req(Wrl)). In fact when $N=2$, both $W_r$ and $W_l$ are isomorphic to $W$ via the following birational correspondences: 
\bea(ll)
W \simeq W_r , & (t, \lambda) \mapsto (t_r, \lambda_r) = (\frac{k t}{1-k'\lambda},  \frac{{\rm i} k \lambda}{1-k'\lambda}), \\
W \simeq W_l , & (t, \lambda) \mapsto (t_l, \lambda_l) = ( \frac{{\rm i}  k t}{k' -\lambda^{-1}}, \frac{1-k'\lambda}{k}), 
\elea(WWrl2)
by which $\varrho_r, \varrho_l$ in (\req(birWrl)) are identified with $\varrho$ in (\req(birWHW)). Since $H$ is a normal subgroup of $\widetilde{G}_2$, $\widetilde{G}_2/H$ gives rise to an automorphism group of ${\goth W}^{(2)}/H$ and $W$. The action of $\widetilde{G}_2$ on $W$ extends the relations in (\req(isHrlH)) through the birational identification (\req(WWrl2)), e.g. the composite $W_{k', k}  \simeq W_{l, k', k} \stackrel{S}{\simeq} W_{k, k'}$ of isomorphisms in (\req(WWrl2)) and (\req(isHrlH)) provides the identification $S:W_{k', k}  \simeq  W_{k, k'}$.
\par \noindent
{\bf Remark.} The birational equivalences in (\req(WWrl2)) hold only in $N=2$. When $N \geq 3$, $W , W_r , W_l$ are not birational equivalent.

\section{Concluding Remarks \label{sec.F}}
In this work, we perform a thorough mathematical investigation of symmetries related to rapidities in CPM within the context of group theory. The set-up is conceptually based upon the analysis of common features of rapidity automorphisms (\req(Sym)) in $N$-state CPM for all $N$, then generalizes the structure to modular symmetries of the rapidity family (\req(WFaml)) that was revealed in the elliptic K3 surface for $N=2$.
By using this approach, the various aspects of the structure groups are studied in Section \ref{sec:Alg}.
Through the representation theory, the structure group can be identified with the automorphism group of all rapidity curves, which constitute the Fermat hypersurface (\req(WFerm)).  In Section \ref{sec.Fermat}, we perform a detailed investigation about the geometrical and symmetry properties of the rapidity fibration of Fermat surface and its associated hyperelliptic-fibered surfaces in CPM, in the context of surface theory in algebraic geometry.  In particular, the intriguing configuration as well as the singularity in degenerated rapidities have developed certain special features in the global geometrical structure of algebraic surfaces involved.  
This paper contains several new observations about the rapidity family of CPM in mathematics and physics, especially those involving the degenerated rapidity curves. Though the structure groups and Fermat rapidity surfaces can be thought of as interesting mathematical topics in their own rights, the physical implications could make it more significant and easier to understand. One relevant physical problem is to find out which ones among $N$-state CPM provide the same theory in statistical mechanics. A pre-condition is the similarity of their rapidities, or equivalently, the rapidity curves are isomorphic under $\widetilde{G}_N$-relations in (\req(GtWf)). Hence the possible temperature candidates are among (\req(24k')). For instance, as hinted at the end of Subsection \ref{subs.hypell}, the modular symmetry $T$  responses the Kramers-Wannier duality of chiral Potts model in \cite{R905}. Along this line, the possibility of other equivalent theories in CPM is currently under investigation.

\end{document}